\documentclass[preprint,12pt]{elsarticle}

\usepackage{amssymb}
\usepackage{amsmath}
\usepackage{amsthm}
\usepackage{eucal}
\usepackage{mathrsfs}
\usepackage{fix-cm}
\usepackage{bm}
\usepackage{amsfonts}
\usepackage{enumerate}
\usepackage{color}

\newcommand{\red}[1]{\textcolor{black}{#1}}

\def\dO{\ \text{d}\hspace{-0.3mm} \Omega}
\def\dG{\ \text{d}\hspace{-0.3mm}\Gamma}

\def\x{\bm{x}}

\def \DD {\mathbb D}

\def \P {\mathscr{P}}

\def \S {\mathscr{S}}
\def \U {\mathscr{U}}

\usepackage{subcaption}
\usepackage{here}

\begin{document}

\begin{frontmatter}
\title{Topology Optimization considering Shielding and Penetrating Features based on Fictitious Physical Model}
\author[inst1]{Daiki Soma}
\author[inst1]{Kota Sakai}
\author[inst1,inst2]{Takayuki Yamada \corref{mycorrespondingauthor}}
\affiliation[inst1]{organization={Department of Mechanical Engineering, Graduate School of Engineering, The University of Tokyo},
            addressline={Yayoi 2--11--16}, 
            city={Bunkyo--ku},
            postcode={113--8656}, 
            state={Tokyo},
            country={Japan}}
\affiliation[inst2]{organization={Department of Strategic Studies, Institute of Engineering Innovation, Graduate School of Engineering, The University of Tokyo},
            addressline={Yayoi 2--11--16}, 
            city={Bunkyo--ku},
            postcode={113--8656}, 
            state={Tokyo},
            country={Japan}}

\cortext[mycorrespondingauthor]{Corresponding author}
\ead{t.yamada@mech.t.u-tokyo.ac.jp}           
\begin{abstract}
This paper proposes topology optimization for considering shielding and penetrating features.
Based on the fictitious physical model, which is a useful approach to control geometric features, the proposed method analyzes fictitious steady-state temperature fields and interprets target geometric features by examining the temperature change.
First, the concept of topology optimization based on the level set method is introduced.
Next, the basic idea of the fictitious physical model for considering geometric features is explained.
Then, the differences between the shielding and penetrating features are clarified, and the fictitious physical model for evaluating these features is proposed.
Furthermore, topology optimization for the minimum mean compliance problem with geometric conditions is formulated.
Finally, 2D and 3D numerical examples are presented to validate the proposed method. 
\end{abstract}
\end{frontmatter}
\section{Introduction}
\label{sec:intro}
Structural optimization is a method to maximize performance for specific objectives by optimizing the configuration of structures from mathematical and mechanical perspectives.
Among structural optimization, topology optimization is the most flexible method, allowing not only changes in the external shape but also modifications in the form, such as the creation or elimination of voids.
Due to its high degree of design freedom, topology optimization has gained attention as an effective method for generating high-performance designs that may surpass designers' imagination.
Consequently, since Bends{\o}e and Kikuchi \cite{bendsoe1988generating} established the engineering applications of topology optimization, it has been used for a wide range of design problems such as mechanical problems \cite{sigmund1997on,holmberg2013stress}, thermal problems \cite{yamada2011a,nakagawa2023level,ogawa2022new,ogawa2022topology}, fluid problems \cite{dilgen2018topology, guan2024topology}, and electromagnetic wave problems \cite{diaz2010a,deng2016topology,murai2023multiscale}.

\red{However}, in practical applications, it is crucial to consider not only the specific physical performance, but also other factors like production costs or manufacturability.
Therefore, topology optimization that incorporates geometric conditions has become an important research area in recent years.
For example, geometric conditions include controlling the minimum \cite{zhang2016minimum,li2023an} and maximum \cite{fernandez2019an,zobaer2023maximum} thickness of a structure to enhance manufacturability. Furthermore, conditions such as suppressing overhangs to enable additive manufacturing without support structures \cite{langelaar2016topology,allaire2017structural} and eliminating undercuts to allow parts to be easily removed from molds in casting processes \cite{zhou2002,xia2010,tajima2024topology} are also included.

\red{As} one of various geometric conditions, this paper introduces a novel topology optimization that considers shielding and penetrating features.
In this paper, ``shielding" refers to a state in which a structure blocks a connection of the void domain between the boundaries, whereas ``penetrating" refers to a state in which a structure allows such a connection.

\red{
Shielding and penetrating features serve various roles.
Regarding the shielding feature, one of its key functions is protecting internal components from external influences.
For example, in electronic devices, external factors such as dust and electromagnetic interference can pose significant challenges.
Dust intrusion can cause short circuits and poor  electrical contact. Moreover, dust absorbs moisture, which may cause corrosion of metal parts and shorten the lifespan of electronic devices.
Similarly, electromagnetic waves can induce noise in circuits, leading to malfunctions and performance degradation.
To prevent these adverse effects, the casing of electronic devices needs to be designed with the shielding feature.
Regarding the penetrating feature, one of its key functions is improving ventilation and facilitating heat dissipation.
For instance, in computers, heat removal is vital for enhancing performance.
To efficiently dissipate the heat generated inside, airflow is essential, and penetrating holes that allow airflow between the inside and outside are required.
Therefore, the casing of computers needs to be designed with the penetrating feature.
As presented in these examples, shielding and penetrating features play important roles, and highly functional designs can be realized by incorporating each feature.}
From this, we propose evaluation models for shielding and penetrating features and apply them to topology optimization using the level set method as a typical method.
\red{It is worth mentioning that this paper is the first implementation to consider shielding and penetrating features in topology optimization.}

\red{The} proposed evaluation models are developed based on the fictitious physical model.
The fictitious physical model is an approach for considering geometric features and introduces fictitious physical fields and governing equations, which are separate from the physical model that represents mechanical phenomena.
Since the governing equations and evaluation functions of this fictitious physical model are formulated in the same way as those of the standard physical model, it enables formulation within the framework of topology optimization.

\red{This} method has gained attention as an approach for extracting geometric features \cite{li2018topology,yamada2019geometric,tajima2023}, and previous research includes the following.
Huang et al. \cite{huang2024}  introduced a fictitious unsteady temperature field to control the length scale of the structures and obtained skeletons with sensitivity information by focusing on the divergence of the temperatures.
Dong et al. \cite{dong2024porous} considered a non-uniform thermal conduction process using skeletons as heat sources and determined the shape of porous structures by projecting the temperature field controlled by the heat coefficient.
Wang et al. \cite{wang2024} incorporated a fictitious heat conduction process and regularized the pseudo-time field to derive a feasible manufacturing sequence for additive manufacturing.

\red{As} described in these examples, the fictitious physical model is useful as an approach to control geometric features. 
Therefore, this paper introduces fictitious steady-state temperature fields as the fictitious physical model, and evaluates shielding and penetrating features by focusing on the temperature gradient.

\red{The} rest of this paper is organized as follows: Section 2 provides an overview of topology optimization based on the level set method. Section 3 explains the idea of the fictitious physical model used to represent geometric features. Section 4 defines the shielding and penetrating features and then proposes evaluation models for these features using the fictitious physical model. Section 5 discusses the minimum mean compliance problem as a concrete example of optimization problems and formulates it including shielding and penetrating conditions. In Section 6, several numerical examples are presented to investigate the effect of the parameters and the usefulness of the proposed method. Finally, Section 7 provides the conclusions of this paper.

\color{black}
\section{Conventional Approaches Adopted in This Paper}
In this paper, the level set method is employed for topology optimization.
Additionally, the fictitious physical model is used to evaluate the geometric features.
The fundamental concepts of these methods are described as follows.
\subsection{Topology Optimization Based on Level Set Method}
\color{black}
In this section, we discuss the key aspects of topology optimization based on the level set method. Topology optimization treats structural optimization problems as solid distribution problems within the fixed design domain, denoted as $D$ \cite{bendsoe1988generating}.
Here, the fixed design domain $D$ is distinguished between the solid domain $\Omega\subset D$, which represents the domain of interest for the design, and the complementary domain $D\setminus \Omega$.
This division is defined using a characteristic function $\chi \in L^{\infty}(D)$ as described below.
\begin{equation}
    \chi(\x) := \begin{cases}
        1 \qquad &{\rm for} \quad \x \in \Omega\\
        0 \qquad &{\rm for} \quad \x \in D \setminus \Omega
    \end{cases}
\end{equation}
This characteristic function is capable of representing structures with arbitrary topology, even allowing infinitesimal structures. 
Therefore, the topology optimization problem is known as an ill-posed problem \cite{allaire2012shape} and must be replaced by a well-posed problem to obtain the optimal solution.
In previous research, the homogenized design method \cite{bendsoe1988generating,suzuki1991homogenization} and the density-based method \cite{bendsoe1999material} are overcome the problem by using relaxation or regularization methods.
In addition to these, the level set-based structural optimization \cite{sethian2000structural,wang2003level,allaire2004structural,luo2008level,wang2018velocity} is also often used, as it provides clear boundaries.
In these methods, the structural shape is represented using the the iso-surface of a distributed scalar function called level set function.
The level set function is defined as follows:
\begin{equation}
    \begin{cases}
        0 < \phi(\x) \le 1 \qquad &{\rm for} \quad \x \in \Omega \setminus \partial \Omega\\
        \phi(\x) = 0 \qquad &{\rm for} \quad \x \in \partial \Omega\\
        -1 \le \phi(\x) < 0 \qquad &{\rm for} \quad \x \in D \setminus \Omega
    \end{cases}
\end{equation}
The level set function $\phi$ takes a positive value in the solid domain, a negative value in the void domain, and the zero iso-surface on the boundary between the two domains.
The characteristic function is as follows using the level set function $\phi$.
\begin{equation}
    \chi_\phi = \begin{cases}
        1 \qquad &{\rm for} \quad \phi(\x) \ge 0\\
        0 \qquad &{\rm for} \quad \phi(\x) < 0
    \end{cases}
\end{equation}
In numerical implementation, the distribution of $\chi_\phi$ needs to smoothly transition between the solid and void domains. To achieve this,  $\chi_\phi$ is re-defined using an approximate Heaviside function as follows:
\begin{equation} \label{eq:chi reg}
    \chi_\phi = \begin{cases}
        1 &{\rm for} \quad \phi(\x)>\delta \\
        \frac{1}{2}+\frac{15}{16} \frac{\phi(\x)}{\delta}-\frac{5}{8}\left(\frac{\phi(\x)}{\delta}\right)^3+\frac{3}{16}\left(\frac{\phi(\x)}{\delta}\right)^5 &{\rm for} \quad \phi(\x) \in\left[-\delta, \delta\right] \\
        0 &{\rm for} \quad \phi(\x)<-\delta
    \end{cases}
\end{equation}
where $\delta >0$ is a parameter representing the transition width of the approximate Heaviside function.

The level set based topology optimization minimizes (or maximizes) the objective function by updating the level set function $\phi$, since it is difficult to directly determine the optimal distribution of the level set function.
To elaborate on this update method for the level set function, we consider the optimization problem defined as follows, where the objective function $J$ is minimized by optimizing the solid distribution in $D$.
\begin{align}
    \inf_{\chi_\phi}~~~J[\chi_\phi]
    \label{eq: opt_problem_lsf}
\end{align}

Here, we introduce the update method proposed by Yamada et al. \cite{yamada2010topology} that is used in this paper.
First, we interpret the optimization problem described by equation (\ref{eq: opt_problem_lsf}) as a reaction-diffusion problem for the level set function.
Specifically, a fictitious time $t$ is introduced for the temporal evolution, and the level set function $\phi$ is updated by the following reaction-diffusion equation:
\begin{align} \label{eq:reaction-diffusion}
    \frac{\partial \phi}{\partial t} = -K \frac{\displaystyle J' \int_D \dO}{\displaystyle \int_D |J'| \dO}
\end{align}
where $K>0$ is a proportionality coefficient.
Given that the optimization problem is ill-posed, a regularization term is added to the equation above to regularize it as follows:
\begin{align} \label{eq:reaction_diffusion_eq}
    \frac{\partial \phi}{\partial t} = -K \left( \frac{\displaystyle J' \int_D \dO}{\displaystyle \int_D |J'| \dO} - \tau \nabla^2 \phi \right)
\end{align}
where $\tau >0$ is the regularization parameter, and increasing the value of $\tau$ suppresses geometric complexity in the optimal shape.
Denoting the level set function representing the initial structure as $\phi_0(\bm{x})$, the level set function is updated by solving the following system.
\begin{eqnarray}
\left\{
\begin{array}{ll}
\cfrac{\partial \phi}{\partial t} = -K \left( \frac{\displaystyle J' \int_D \dO}{\displaystyle \int_D |J'| \dO} - \tau \nabla^2 \phi \right) & \mathrm{in}~~D \\
\phi=1 & \mathrm{on}~~\partial D_m \\
\bm{n} \cdot \nabla \phi = 0 & \mathrm{on}~~\partial D \setminus \partial D_m \\
\left. \phi(\bm{x}) \right|_{t = 0}= \phi_0(\bm{x})
\end{array}
\right. 
\end{eqnarray}
where $\bm{n}$ is the outward unit normal vector on the boundary $\partial D$, and the boundary $\partial D_m$ is connected to the solid domain outside of the fixed design domain $D$.
\color{black}
\subsection{Fictitious Physical Model for Geometric Feature Evaluation}
\color{black}
In general, topology optimization involves formulating objective functions and governing equations based on the standard physical model.
The standard physical model refers to a set of governing equations that describe actual physical phenomena.
Typical examples include the elasticity equation for solid deformation and the Navier-Stokes equations for fluid flow.
In topology optimization, state variables are derived from numerically solvable governing equations, which are then used to formulate evaluation indices such as the objective function. Hence, while the standard physical model enables obtaining optimal designs from a mechanical perspective, incorporating geometric conditions is challenging. This difficulty arises from the lack of state variables that explicitly  represent geometric features.
Therefore, the fictitious physical model has been proposed to incorporate geometric conditions.
The fictitious physical model is the fictitious system of governing equations that describes geometric features, and formulates evaluation indices for target geometric features by introducing fictitious physical fields.
For example, 
Yamada et al. \cite{yamada2022exclusion} adopted steady-state thermal diffusion equations assuming that the solid domain is adiabatic and the void domain is the heat source in order to exclude closed cavities in structures, which pose manufacturing issues in additive manufacturing. 
Sato et al. \cite{sato2017manufacturability} employed steady-state anisotropic advection-diffusion equations assuming the fictitious light irradiation in the solid domain to obtain structures that allow the removal of molded parts in casting.
Incorporating the fictitious physical model enables the formulation of the following multi-physics problems from both mechanical and geometric perspectives.
\begin{align}
    \inf_{\chi_\phi} \qquad &J[x_1, x_2, \ldots, x_m,\: y_1, y_2, \ldots, y_n,\: \chi_\phi] \notag\\
    \textrm{subject to:} \qquad &\textrm{Governing equations for $x_1, x_2, \ldots, x_m$} \notag\\
                                &\textrm{Governing equations for $y_1,\, y_2, \ldots,\, y_n$} \notag
\end{align}
where $x_1, x_2, \ldots, x_m$ and $y_1,\, y_2, \ldots,\, y_n$ represent the state variables for the standard physical model and the fictitious physical model, respectively. 
Using this approach, it becomes feasible to optimize designs that account for geometric conditions within the basic framework of topology optimization.
\section{Evaluation of Target Geometric Feature}

Shielding and penetrating features are required in various situations in structural design. The shielding feature contributes to protection from the outside and the prevention of leakage. For example, casings of electronic devices shield sensitive internal components from external electromagnetic interference, and cylinder covers help maintain internal pressure, ensuring safe operation.
In contrast, the penetrating feature supports improved  ventilation and spatial integration. For instance, in computer and automotive cooling systems, heat-generating components are placed in the same space as ventilation openings for efficient heat removal. As demonstrated, it is important to consider shielding and penetrating features in various designs.
Therefore, this paper focuses on the geometric features related to shielding and penetrating features.

\begin{figure}[H]
	\begin{center}
		\includegraphics[height=4cm]{ 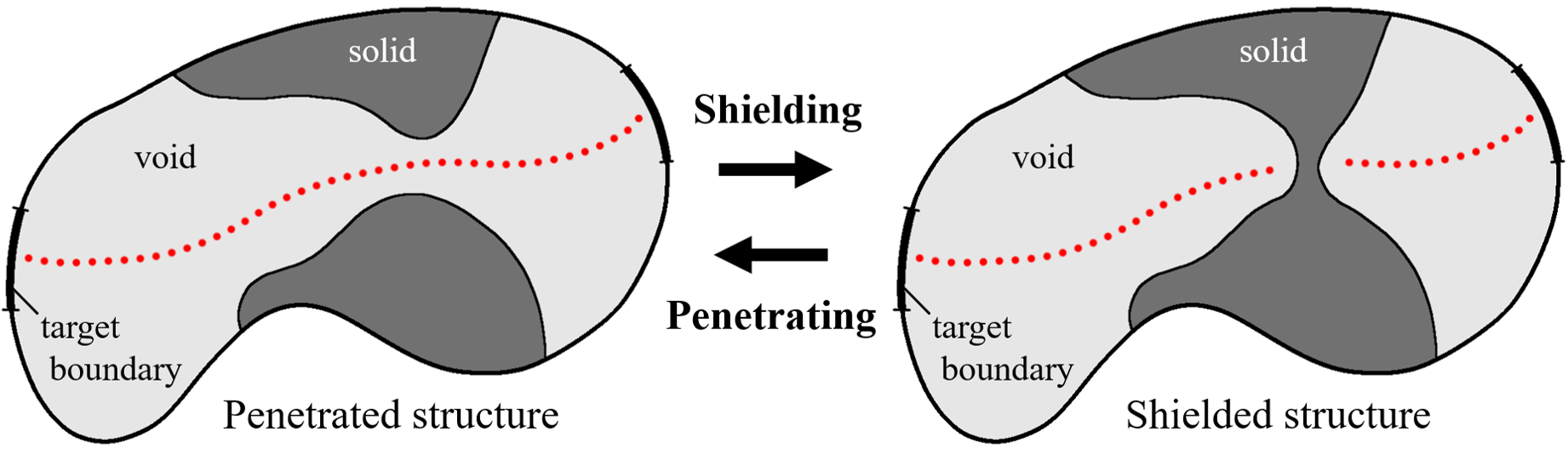}
		\caption{Definition of shielded and penetrated structures}
		\label{fig:definition}
	\end{center}
\end{figure}

The target geometric features are defined as illustrated in Figure \ref{fig:definition}.
In the case of the penetrated structure, the two sets of target boundaries are connected through the same void domain, as shown by the red dotted line, whereas they are not connected in the case of the shielded structure.
To evaluate these features, we formulate them using the fictitious physical model.

\red{
Section~\ref{sec:Boun} examines the boundary conditions suitable for the evaluation model.}
Section~\ref{sec:Evaluation S} presents the evaluation model for the shielding feature.
Section~\ref{sec:Evaluation P} provides the evaluation model for the penetrating feature.


\color{black}
\subsection{Influence of Boundary Conditions}
\label{sec:Boun}
This section investigates boundary conditions suitable for evaluating the shielding and penetrating features.
The concept of this study is to evaluate the target geometric features by focusing on differences in temperature gradients within a steady-state temperature field \(T\).
In the shielding and penetrating evaluation, the following points are required:

\begin{itemize}
    \item
    For the shielding evaluation, it is important to identify the void regions that hinder the shielding feature.
    Therefore, it is necessary to ensure the significant difference between the temperature gradient $T^\mathrm{pene}_\mathrm{void}$ in the void domain of a penetrated structure and the temperature gradient $T^\mathrm{shie}_\mathrm{void}$ in the void domain of a shielded structure.
    \item 
    For the penetrating evaluation, it is important to identify the solid regions that hinder the penetrating feature.
    Therefore, it is necessary to ensure the significant difference between the temperature gradient $T^\mathrm{pene}_\mathrm{solid}$ in the solid domain of a penetrated structure and the temperature gradient $T^\mathrm{shie}_\mathrm{solid}$ in the solid domain of a shielded structure.
\end{itemize}
In order to satisfy these requirements, the boundary conditions must be properly set in the steady-state thermal analysis.

\begin{figure}[H]
    \begin{center}
        \begin{minipage}[t]{0.49\columnwidth}
            \centering
            \includegraphics[width=1\columnwidth]{ 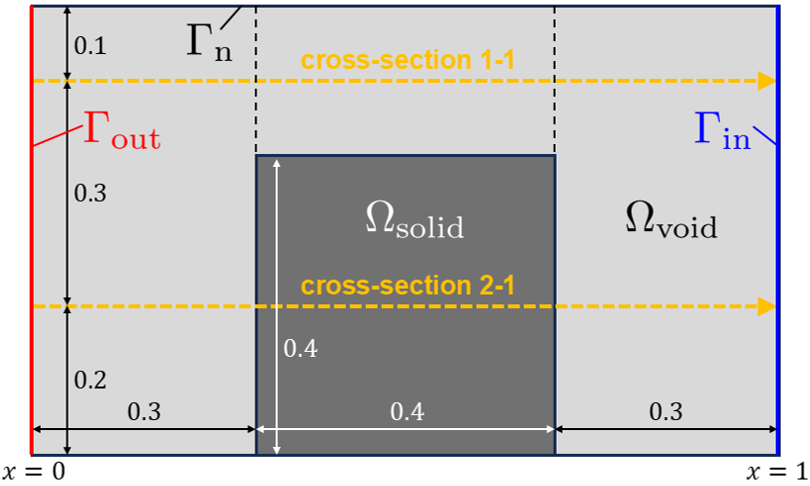}
            \subcaption{Penetrated structure and specific cross-sections}
        \end{minipage}
        \begin{minipage}[t]{0.49\columnwidth}
            \centering
            \includegraphics[width=1\columnwidth]{ 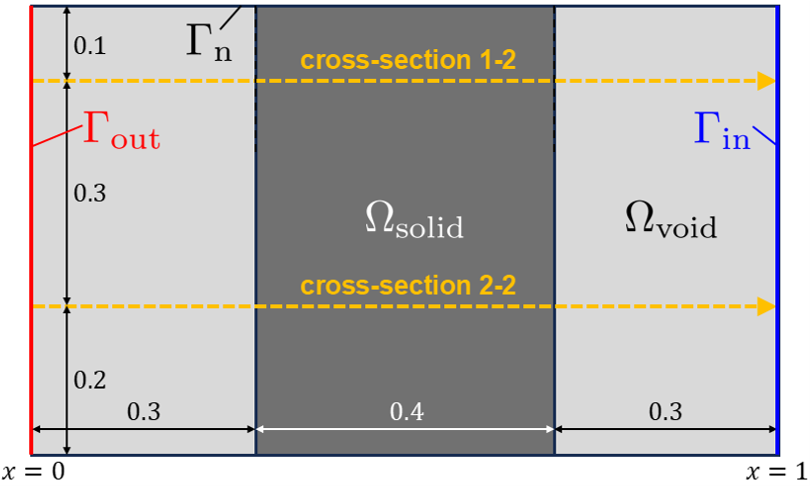}
            \subcaption{Shielded structure and specific cross-sections}
        \end{minipage}
        \caption{Problem setting for thermal analysis.}
	\label{fig:analysis}
    \end{center}
\end{figure}

To examine the influence of boundary conditions, as a specific example, we consider the problem setting shown in Figure \ref{fig:analysis}, which includes a penetrated structure (a) and a shielded structure (b) composed of the solid domain $\Omega_\mathrm{solid}$ and the void domain $\Omega_\mathrm{void}$.
We analyze the temperature distributions on specific cross-sections of Figure \ref{fig:analysis}, based on the following two governing equations.
\begin{align}
    \begin{cases} \label{eq:GovDi}
        -\nabla \cdot \left(a \nabla T \right) = 0 \, \, \, \, & \text{in} \, \, \, \, \Omega_\mathrm{solid} \cup \Omega_\mathrm{void}\\
         T = 1 \, \, \, \, & \text{on} \, \, \, \, \Gamma_\mathrm{out}\\
        T = 0 \, \, \, \, & \text{on} \, \, \, \, \Gamma_\mathrm{in}\\
        -\bm{n} \cdot \left(a \nabla T \right) = 0 \qquad & \text{on}  \, \, \, \, \Gamma_\mathrm{n}
    \end{cases}
\end{align}
\begin{align}
    \begin{cases} \label{eq:GovNo}
        -\nabla \cdot \left(a \nabla T \right) = 0 \, \, \, \, & \text{in} \, \, \, \, \Omega_\mathrm{solid} \cup \Omega_\mathrm{void}\\
         -\bm{n} \cdot \left(a \nabla T \right) = 1 \, \, \, \, & \text{on} \, \, \, \, \Gamma_\mathrm{out}\\
        T = 0 \, \, \, \, & \text{on} \, \, \, \, \Gamma_\mathrm{in}\\
        -\bm{n} \cdot \left(a \nabla T \right) = 0 \qquad & \text{on}  \, \, \, \, \Gamma_\mathrm{n}
    \end{cases}
\end{align}
where \(a\) is a coefficient defined as follows:
\begin{equation}
a = \left\{
\begin{array}{cl}
1 & \mathrm{in}~~\Omega_\mathrm{solid} \\
100 & \mathrm{in}~~\Omega_\mathrm{void}
\end{array}
\right.
\end{equation}
The difference is that the governing equation (\ref{eq:GovDi}) applies the Dirichlet boundary condition to $\Gamma_\mathrm{out}$, while the governing equation (\ref{eq:GovNo}) applies the Neumann boundary condition.

\begin{figure}[H]
    \centering
    \begin{minipage}[t]{0.49\linewidth}
        \centering
        \includegraphics[width=\linewidth]{ 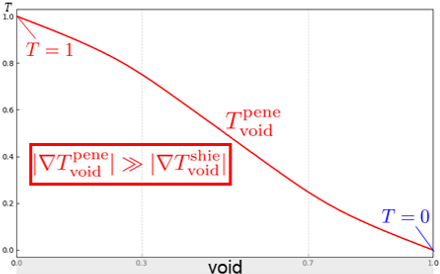}
        \subcaption{Distribution on cross-section 1-1 obtained from governing equation (\ref{eq:GovDi})}
    \end{minipage}
    \begin{minipage}[t]{0.49\linewidth}
        \centering
        \includegraphics[width=\linewidth]{ 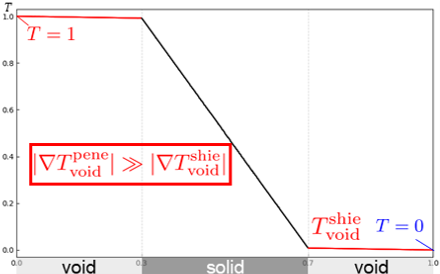}
        \subcaption{Distribution on cross-section 1-2 obtained from governing equation (\ref{eq:GovDi})}
        \vspace{0.02\linewidth}
    \end{minipage}
    \begin{minipage}[t]{0.49\linewidth}
        \centering
        \includegraphics[width=\linewidth]{ 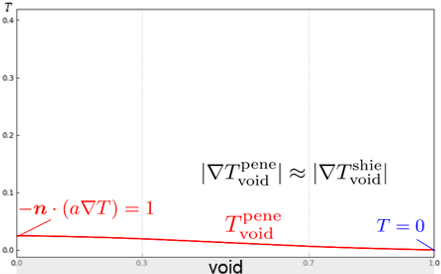}
        \subcaption{Distribution on cross-section 1-1 obtained from governing equation (\ref{eq:GovNo})}
    \end{minipage}
    \begin{minipage}[t]{0.49\linewidth}
        \centering
        \includegraphics[width=\linewidth]{ 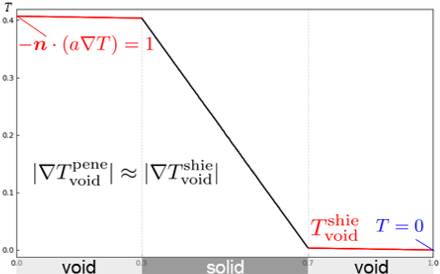}
        \subcaption{Distribution on cross-section 1-2 obtained from governing equation (\ref{eq:GovNo})}
    \end{minipage}
    \caption{Temperature distributions with different boundary conditions (Analysis of the shielding condition)}
    \label{fig:distributionsS}
\end{figure}

To examine the boundary conditions suitable for the shielding evaluation, Figure \ref{fig:distributionsS} shows the temperature distributions on cross-sections 1-1 and 1-2, derived from governing equations (\ref{eq:GovDi}) and (\ref{eq:GovNo}).
Under the boundary conditions of governing equation (\ref{eq:GovNo}), $|T^\mathrm{pene}_\mathrm{void}|$ in (c) and $|T^\mathrm{shie}_\mathrm{void}|$ in (d) are similar. However, under the boundary conditions of governing equation (\ref{eq:GovDi}), $|T^\mathrm{pene}_\mathrm{void}|$ in (a) is noticeably larger than $|T^\mathrm{shie}_\mathrm{void}|$ in (b).
This indicates that the boundary conditions of governing equation (\ref{eq:GovDi}) are appropriate for the shielding evaluation.

\begin{figure}[H]
    \centering
    \begin{minipage}[t]{0.49\linewidth}
        \centering
        \includegraphics[width=\linewidth]{ 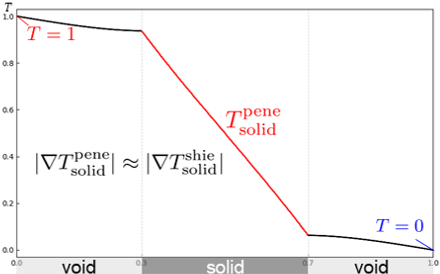}
        \subcaption{Distribution on cross-section 2-1 obtained from governing equation (\ref{eq:GovDi})}
    \end{minipage}
    \begin{minipage}[t]{0.49\linewidth}
        \centering
        \includegraphics[width=\linewidth]{ 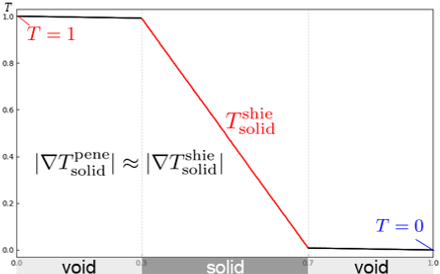}
        \subcaption{Distribution on cross-section 2-2 obtained from governing equation (\ref{eq:GovDi})}
        \vspace{0.02\linewidth}
    \end{minipage}
    \begin{minipage}[t]{0.49\linewidth}
        \centering
        \includegraphics[width=\linewidth]{ 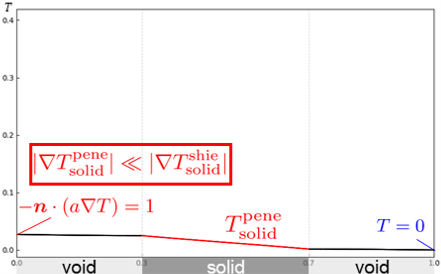}
        \subcaption{Distribution on cross-section 2-1 obtained from governing equation (\ref{eq:GovNo})}
    \end{minipage}
    \begin{minipage}[t]{0.49\linewidth}
        \centering
        \includegraphics[width=\linewidth]{ 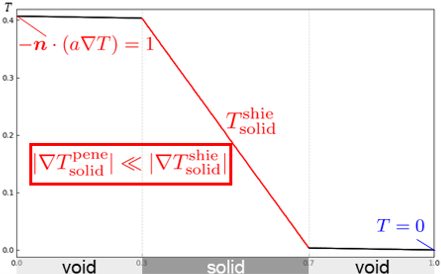}
        \subcaption{Distribution on cross-section 2-2 obtained from governing equation (\ref{eq:GovNo})}
    \end{minipage}
    \caption{Temperature distributions with different boundary conditions (Analysis of the penetrating condition)}
    \label{fig:distributionsP}
\end{figure}

To investigate the boundary conditions suitable for the penetrating evaluation, Figure \ref{fig:distributionsP} presents the temperature distributions on cross-section 2-1 and 2-2, derived from governing equations (\ref{eq:GovDi}) and (\ref{eq:GovNo}).
Under the boundary conditions of governing equation (\ref{eq:GovDi}), $|T^\mathrm{shie}_\mathrm{solid}|$ in (b) and $|T^\mathrm{pene}_\mathrm{solid}|$ in (a) are similar. However, under the boundary conditions of governing equation (\ref{eq:GovDi}), $|T^\mathrm{shie}_\mathrm{solid}|$ in (d) is significantly larger than $|T^\mathrm{pene}_\mathrm{solid}|$ in (c).
This suggests that the boundary conditions of the governing equation (\ref{eq:GovNo}) are appropriate for the penetrating evaluation.

Based on the above, in this paper, the Dirichlet condition is applied as the heat source for the shielding evaluation, and the Neumann condition is applied as the heat source for the penetrating evaluation.

\color{black}
\subsection{Evaluation of Shielding Feature} \label{sec:Evaluation S}
To evaluate the shielding feature, we introduce the fictitious physical variable $s \in H^{1}(D)$ and define its governing equation as follows.
\begin{align}
    \begin{cases} \label{eq:GovS}
        -\nabla \cdot \left(a_s \nabla s \right) = 0 \, \, \, \, & \text{in} \, \, \, \, D\\
        s = 1 \, \, \, \, & \text{on} \, \, \, \, \Gamma^{\,\mathrm{s}}_\mathrm{out}\\
        s = 0 \, \, \, \, & \text{on} \, \, \, \, \Gamma^{\,\mathrm{s}}_\mathrm{in}\\
        -\bm{n} \cdot \left(a_s \nabla s \right) = 0 \qquad & \text{on}  \, \, \, \, \partial D \setminus \left(\Gamma^{\,\mathrm{s}}_\mathrm{out} \cup \Gamma^{\,\mathrm{s}}_\mathrm{in} \right)
    \end{cases}
\end{align}
where
\begin{align}
    a_s := \left\{ \kappa^{\,\mathrm{s}}_\mathrm{solid} \chi_\phi^3 + \kappa^{\,\mathrm{s}}_\mathrm{void} \left(1 - \chi_\phi^3 \right) \right\} L^2_s
\end{align}
\red{
is a non-dimensionalized diffusion coefficient.
This non-dimensionalization allows consistent diffusion coefficients to be applied regardless of the design domain's scale.}
$\kappa^{\,\mathrm{s}}_\mathrm{solid}$ and $\kappa^{\,\mathrm{s}}_\mathrm{void}$ represent the diffusion coefficients in the solid and void domains, respectively. These coefficients are set such that $0 < \kappa^{\,\mathrm{s}}_\mathrm{solid} < \kappa^{\,\mathrm{s}}_\mathrm{void}$. The ratio of the coefficients $\kappa^{\,\mathrm{s}}_\mathrm{solid}$ and  $\kappa^{\,\mathrm{s}}_\mathrm{void}$ is examined in the numerical examples in Section~\ref{sec:paramS effect}. The characteristic length $L_s$ is defined as the typical length between the boundaries 
$\Gamma^{\,\mathrm{s}}_\mathrm{out}$ and $\Gamma^{\,\mathrm{s}}_\mathrm{in}$.
\red{
In addition, the exponent of $\chi_\phi$ is set to 3 in order to utilize the SIMP method framework in the sensitivity analysis.
}

\begin{figure}[H]
    \centering
    \begin{minipage}[t]{0.46\linewidth}
        \centering
        \includegraphics[width=\linewidth]{ 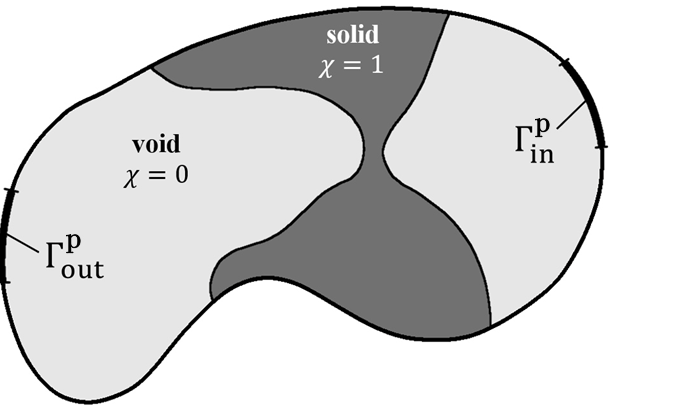}
        \subcaption{Problem setting regarding penetrated structure}
    \end{minipage}
    \hspace{0.03\linewidth}
    \begin{minipage}[t]{0.46\linewidth}
        \centering
        \includegraphics[width=\linewidth]{ 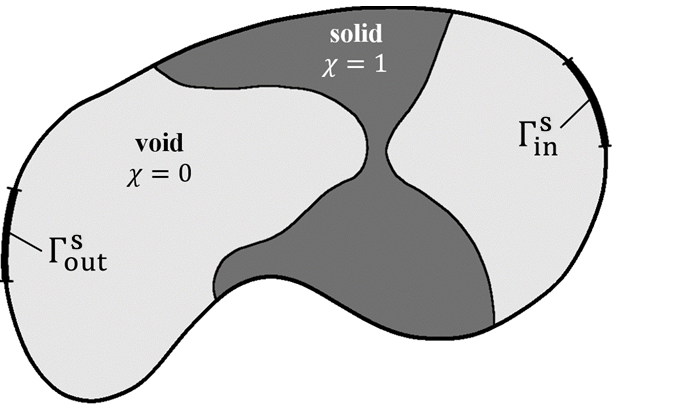}
        \subcaption{Problem setting regarding shielded structure}
    \end{minipage}
    \vspace{0.02\linewidth}
    \begin{minipage}[t]{0.46\linewidth}
        \centering
        \includegraphics[width=\linewidth]{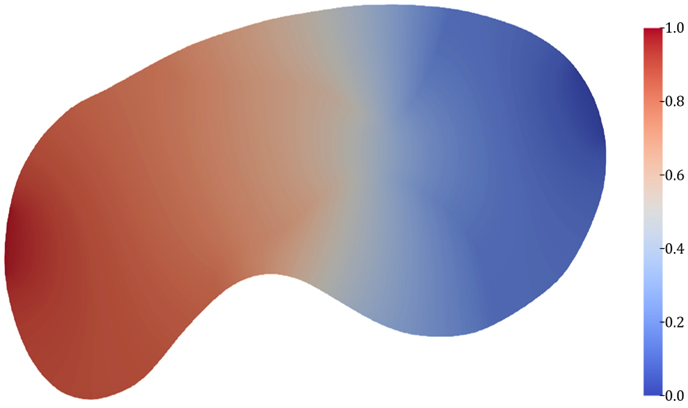}
        \subcaption{Distribution of $s$ in the case of penetrated structure}
    \end{minipage}
    \hspace{0.03\linewidth}
    \begin{minipage}[t]{0.46\linewidth}
        \centering
        \includegraphics[width=\linewidth]{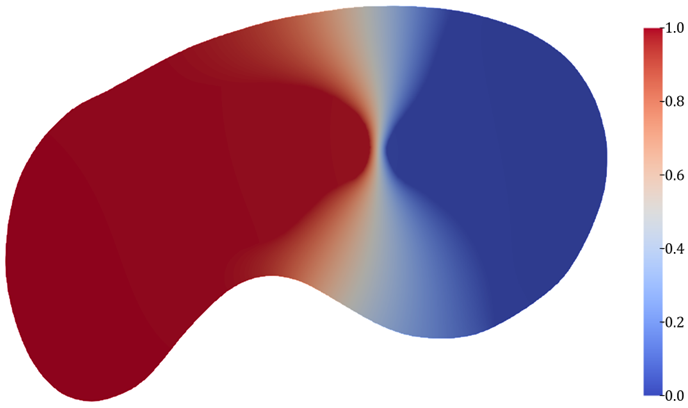}
        \subcaption{Distribution of $s$ in the case of shielded structure}
    \end{minipage}
    \vspace{0.02\linewidth}
    \begin{minipage}[t]{0.46\linewidth}
        \centering
        \includegraphics[width=\linewidth]{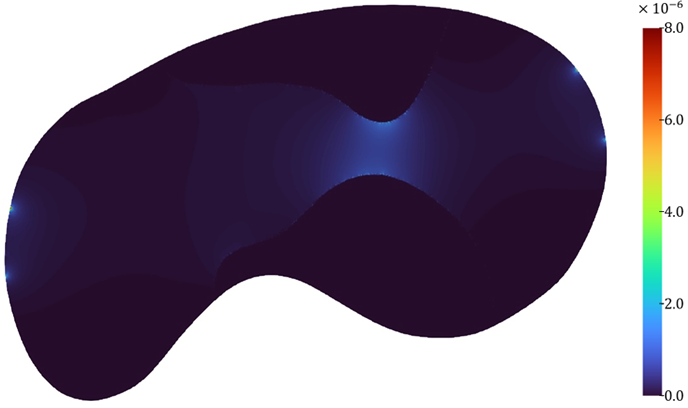}
        \subcaption{Distribution of $(1-\chi_\phi)^3 \nabla s \cdot \nabla s$ in the case of penetrated structure}
    \end{minipage}
    \hspace{0.02\linewidth}
    \begin{minipage}[t]{0.46\linewidth}
        \centering
        \includegraphics[width=\linewidth]{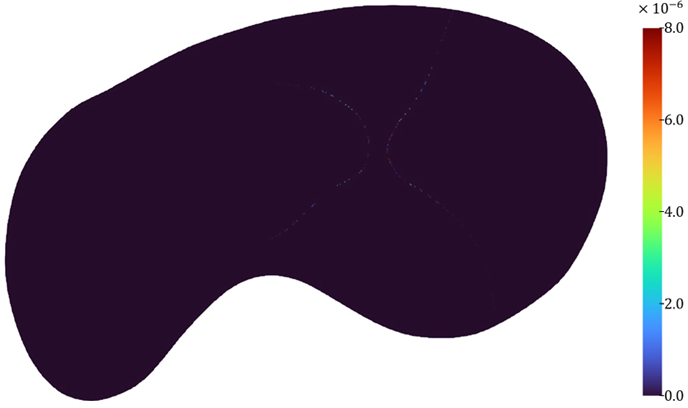}
        \subcaption{Distribution of $(1-\chi_\phi)^3 \nabla s \cdot \nabla s$ in the case of shielded structure}
    \end{minipage}
    \caption{Comparison of distribution between penetrated and shielded structures in terms of shielding feature.}
    \label{fig:evaluS}
\end{figure}

The governing equation (\ref{eq:GovS}) represents steady-state thermal diffusion when a fixed temperature difference is applied between the boundaries $\Gamma^{\,\mathrm{s}}_\mathrm{out}$ and $\Gamma^{\,\mathrm{s}}_\mathrm{in}$.
This model sets the high diffusion coefficient $\kappa^{\,\mathrm{s}}_\mathrm{void}$ in the void domain, which makes the void domain less subject to temperature changes than the solid domain.
Consequently, when the solid domain shields between the boundaries $\Gamma^{\,\mathrm{s}}_\mathrm{out}$ and $\Gamma^{\,\mathrm{s}}_\mathrm{in}$, the temperature change in the solid domain becomes significant, and the temperature change in the void domain becomes nearly zero.
On the other hand, when the solid domain does not shield between the boundaries, a fixed temperature difference arises in the void domain, resulting in a larger temperature change in the unshielded void domain.
To observe the differences in the behavior of $s$ resulting from variations in shielding features, we refer to Figure \ref{fig:evaluS}.
(c) and (d) show the distributions of $s$ derived from the governing equation (\ref{eq:GovS}) for the penetrated structure (a) and the shielded structure (b), respectively, with diffusion coefficients set as $\kappa^{\,\mathrm{s}}_\mathrm{solid} = 1$ and $\kappa^{\,\mathrm{s}}_\mathrm{void} = 1000$.
As stated, setting the diffusion coefficients such that $0 < \kappa^{\,\mathrm{s}}_\mathrm{solid} < \kappa^{\,\mathrm{s}}_\mathrm{void}$ makes the temperature change smaller in the void domain compared to the solid domain.
This phenomenon is particularly prominent in the shielded structure, where the temperature change in the void domain is nearly zero. In contrast, the penetrated structure exhibits large temperature change in the void domain.
From the different behaviors of $s$, we define the following evaluation function $J_s$, focusing on temperature change in the void domain as an indicator of the shielding feature.
\begin{equation} \label{eq:General Js}
	J_s:= \int_D \left(1-\chi_\phi\right)^{c_s} f_s(|\nabla s|) \dO
\end{equation}
where $f_s$ is a monotonically increasing function, which is set as $f_s(|\nabla s|)=\nabla s \cdot \nabla s$ in this paper. Additionally, $c_s$ is a constant greater than 1 to retain the information of the factor $(1-\chi_\phi)$ as the sensitivity of $J_s$, which is set as $c_s = 3$ in this paper. 
Therefore, this paper expresses the evaluation function $J_s$ as follows:
\begin{equation} \label{eq:Js}
	J_s = \int_D \left(1-\chi_\phi\right)^3 \nabla s \cdot \nabla s \dO
\end{equation}
This evaluation function assesses the shielding feature by examining the temperature gradient in the void domain. 
When a structure has the shielding feature, the temperature gradient in the void domain is nearly zero. Conversely, when a structure does not have the shielding feature, the temperature gradient in the unshielded void domain increases. 
As a result, this evaluation function takes a value close to zero when the shielding is achieved and a large value when the shielding is not achieved.
Verification examples of this function are shown in Figure \ref{fig:evaluS}.
(e) and (f) are the distributions of $\left(1-\chi_\phi\right)^3 \nabla s \cdot \nabla s$ for the penetrated structure (a) and the shielded structure (b), respectively, where the diffusion coefficients are $\kappa^{\,\mathrm{s}}_\mathrm{solid} = 1$ and $\kappa^{\,\mathrm{s}}_\mathrm{void} = 1000$.
It can be seen that the value of $\left(1-\chi_\phi\right)^3 \nabla s \cdot \nabla s$ is larger for the penetrated structure compared to the shielded structure. Especially, the value is larger in the narrower void domain.
Therefore, the shielding effect on the design objects is enhanced by incorporating the evaluation function $J_s$ into the objective function of the minimization problem.

The sensitivity $J_s'$ of the evaluation function (\ref{eq:Js}) is formulated as follows:
\begin{equation} \label{eq:Js'}
J_s'= 3 \left(1-\chi_\phi\right)^2 \nabla s \cdot \nabla s - \frac{\partial a_s}{\partial \chi_\phi} \nabla s \cdot \nabla {\hat{s} }
\end{equation}
where the adjoint variable $\hat{s} \in H^{1}(D)$ is defined as follows:
\begin{align} \label{eq:AdjS}
    \begin{cases}
        -\nabla \cdot (a_s \nabla \hat{s}) = 2 (1 - \chi_\phi)^3 \nabla ^2 s \, \, \, \, & \text{in} \, \, \, \, D\\
        \hat{s} = 0 \, \, \, \, & \text{on} \, \, \, \, \Gamma^{\,\mathrm{s}}_\mathrm{out} \cup \Gamma^{\,\mathrm{s}}_\mathrm{in}\\
        -\bm{n} \cdot \left(a_s \nabla \hat{s} \right) = 0 \qquad & \text{on}  \, \, \, \, \partial D \setminus \left(\Gamma^{\,\mathrm{s}}_\mathrm{out} \cup \Gamma^{\,\mathrm{s}}_\mathrm{in} \right)
    \end{cases}
\end{align}
\subsection{Evaluation of Penetrating Feature} \label{sec:Evaluation P}
To evaluate the penetrating feature, we introduce the fictitious physical variable $p \in H^{1}(D)$ and define its governing equation as follows.
\begin{align} \label{eq:GovP}
    \begin{cases}
        -\nabla \cdot \left(a_p \nabla p \right) = 0 \, \, \, \, & \text{in} \, \, \, \, D\\
        -\bm{n} \cdot \left(a_p \nabla p \right) = 1 \qquad & \text{on} \, \, \, \, \Gamma^{\,\mathrm{p}}_\mathrm{out}\\
        p = 0 \, \, \, \, & \text{on} \, \, \, \, \Gamma^{\,\mathrm{p}}_\mathrm{in}\\
        -\bm{n} \cdot \left(a_p \nabla p \right) = 0 \qquad & \text{on}  \, \, \, \, \partial D \setminus \left(\Gamma^{\,\mathrm{p}}_\mathrm{out} \cup \Gamma^{\,\mathrm{p}}_\mathrm{in} \right)
    \end{cases}
\end{align}
where
\begin{align}
    a_p := \left\{ \kappa^{\,\mathrm{p}}_\mathrm{solid} \chi_\phi^3 + \kappa^{\,\mathrm{p}}_\mathrm{void} \left(1 - \chi_\phi^3 \right) \right\} L^2_p
\end{align}
\red{
is a non-dimensionalized diffusion coefficient.
This non-dimensionalization allows consistent diffusion coefficients to be applied regardless of the design domain's scale.}
$\kappa^{\,\mathrm{p}}_\mathrm{solid}$ and $\kappa^{\,\mathrm{p}}_\mathrm{void}$ represent the diffusion coefficients in the solid and void domains, respectively. These coefficients are set such that $0 < \kappa^{\,\mathrm{p}}_\mathrm{solid} < \kappa^{\,\mathrm{p}}_\mathrm{void}$. The ratio of the coefficients $\kappa^{\,\mathrm{p}}_\mathrm{solid}$ and  $\kappa^{\,\mathrm{p}}_\mathrm{void}$ is examined in the numerical examples in Section~\ref{sec:paramP effect}. The characteristic length $L_p$ is defined as the typical length between the boundaries $\Gamma^{\,\mathrm{p}}_\mathrm{out}$ and $\Gamma^{\,\mathrm{p}}_\mathrm{in}$.
\red{
In addition, the exponent of $\chi_\phi$ is set to 3 in order to utilize the SIMP method framework in the sensitivity analysis.
}

\begin{figure}[H]
    \centering
    \begin{minipage}[t]{0.46\linewidth}
        \centering
        \includegraphics[width=\linewidth]{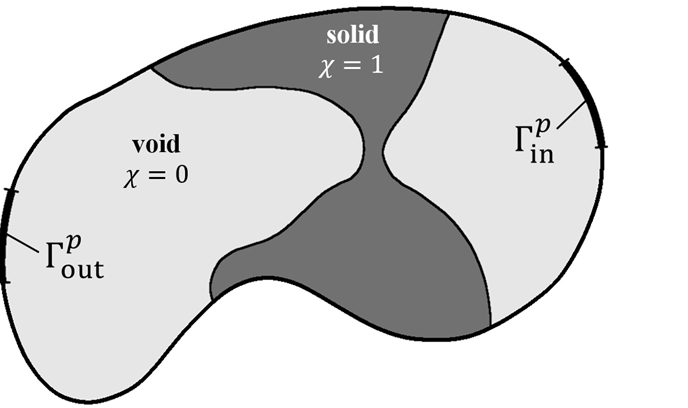}
        \subcaption{Problem setting regarding shielded structure}
    \end{minipage}
    \hspace{0.03\linewidth}
    \begin{minipage}[t]{0.46\linewidth}
        \centering
        \includegraphics[width=\linewidth]{ 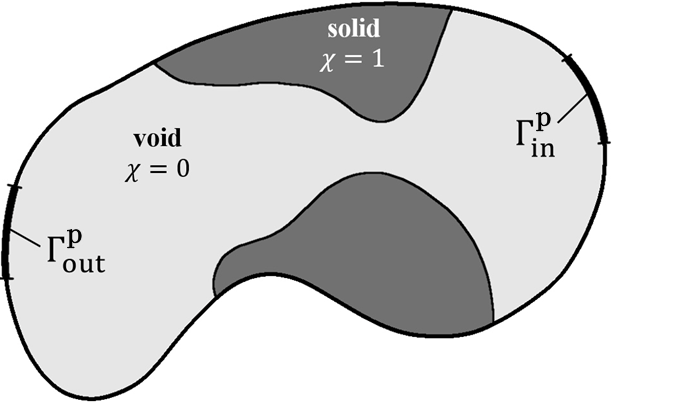}
        \subcaption{Problem setting regarding penetrated structure}
    \end{minipage}
    \vspace{0.02\linewidth}
    \begin{minipage}[t]{0.46\linewidth}
        \centering
        \includegraphics[width=\linewidth]{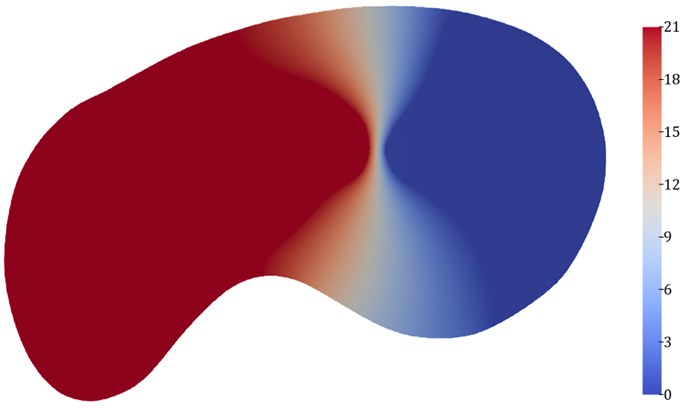}
        \subcaption{Distribution of $p$ in the case of shielded structure}
    \end{minipage}
    \hspace{0.03\linewidth}
    \begin{minipage}[t]{0.46\linewidth}
        \centering
        \includegraphics[width=\linewidth]{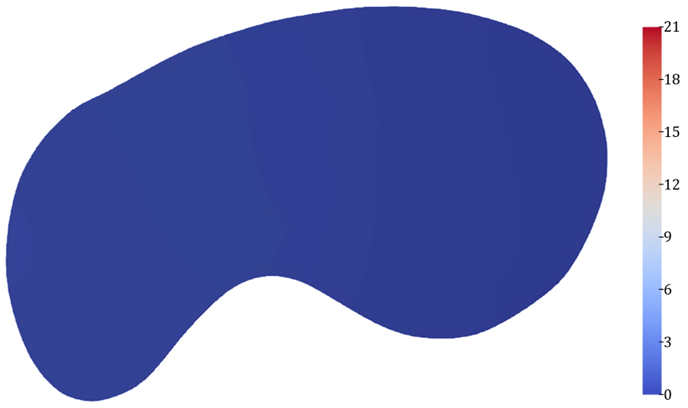}
        \subcaption{Distribution of $p$ in the case of penetrated structure}
    \end{minipage}
    \vspace{0.02\linewidth}
    \begin{minipage}[t]{0.46\linewidth}
        \centering
        \includegraphics[width=\linewidth]{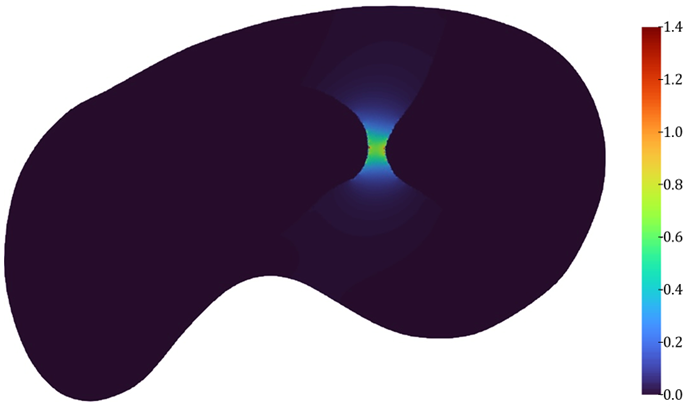}
        \subcaption{Distribution of $\chi_\phi^3 \nabla p \cdot \nabla p$ in the case of shielded structure}
    \end{minipage}
    \hspace{0.02\linewidth}
    \begin{minipage}[t]{0.46\linewidth}
        \centering
        \includegraphics[width=\linewidth]{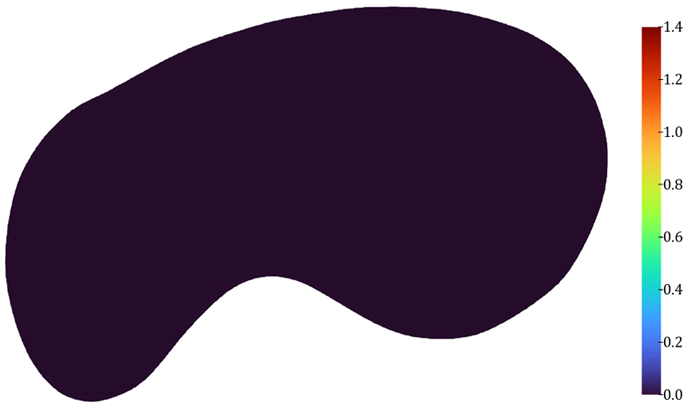}
        \subcaption{Distribution of $\chi_\phi^3 \nabla p \cdot \nabla p$ in the case of penetrated structure}
    \end{minipage}
    \caption{Comparison of distribution between shielded and penetrated structures in terms of penetrating feature.}
    \label{fig:evaluP}
\end{figure}

The governing equation (\ref{eq:GovP}) represents steady-state thermal diffusion when a constant heat flux is applied from the boundary $\Gamma^{\,\mathrm{p}}_\mathrm{out}$ to the boundary $\Gamma^{\,\mathrm{p}}_\mathrm{in}$.
This model sets the high diffusion coefficient $\kappa^{\,\mathrm{p}}_\mathrm{void}$ in the void domain, similar to the model for shielding, which makes the void domain less subject to temperature changes than the solid domain.
Consequently, when the void domain penetrates between the boundaries $\Gamma^{\,\mathrm{p}}_\mathrm{out}$ and $\Gamma^{\,\mathrm{p}}_\mathrm{in}$, heat supplied at the boundary $\Gamma^{\,\mathrm{p}}_\mathrm{out}$
quickly transfers through the void domain to the boundary $\Gamma^{\,\mathrm{p}}_\mathrm{in}$, leading to a small overall temperature change.
As a result, the temperature change in the solid domain is nearly zero.
On the other hand, when the void domain does not penetrate between the boundaries, heat tends to accumulate in the solid domain, resulting in a larger temperature change in the unpenetrated solid domain.
The reason for this behavior is that the temperature difference between the boundaries $\Gamma^{\,\mathrm{p}}_\mathrm{out}$ and $\Gamma^{\,\mathrm{p}}_\mathrm{in}$ is not fixed in the governing equation (\ref{eq:GovP}) for penetrating. This setting is a major difference from the governing equation (\ref{eq:GovS}) for shielding.
To observe the differences in the behavior of $p$ resulting from variations in penetrating features, we refer to Figure \ref{fig:evaluP}.
(c) and (d) show the distributions of $p$ derived from the governing equation (\ref{eq:GovP}) for the shielded structure (a) and the penetrated structure (b), respectively, with diffusion coefficients set as $\kappa^{\,\mathrm{p}}_\mathrm{solid}=1$ and $\kappa^{\,\mathrm{p}}_\mathrm{void}=1000$. 
As stated, by setting the diffusion coefficients such that $0 < \kappa^{\,\mathrm{p}}_\mathrm{solid} < \kappa^{\,\mathrm{p}}_\mathrm{void}$, the void domain can transfer heat more quickly than the solid domain.
Therefore, the temperature is almost uniform overall in the penetrated structure, where the temperature change in the solid domain is nearly zero.
In contrast, the shielded structure exhibits large temperature change in the solid domain.
From the different behaviors of $p$, we define the following evaluation function $J_p$, focusing on the temperature change in the solid domain as an indicator of the penetrating feature.
\begin{equation} \label{eq:General Jp}
	J_p:= \int_D \chi_\phi^{c_p} f_p(|\nabla p|) \dO
\end{equation}
where $f_p$ is a monotonically increasing function, which is set as $f_p(|\nabla p|)=\nabla p \cdot \nabla p$ in this paper. Additionally, $c_p$ is a constant greater than 1 to retain the information of the factor $\chi_\phi$ as the sensitivity of $J_p$, which is set as $c_p = 3$ in this paper. 
Therefore, this paper expresses the evaluation function $J_p$ as follows:
\begin{equation} \label{eq:Jp}
	J_p = \int_D \chi_\phi^3 \nabla p \cdot \nabla p \dO
\end{equation}
This evaluation function assesses the penetrating feature by examining the temperature gradient in the solid domain. 
When a structure has the penetrating feature, the temperature gradient in the solid domain is nearly zero. Conversely, when a structure does not have the penetrating feature, the temperature gradient in the  unpenetrated solid domain increases. 
As a result, this evaluation function takes a value close to zero when the penetrating is achieved and a large value when the penetrating is not achieved.
Verification examples of this function are shown in Figure \ref{fig:evaluP}.
(e) and (f) are the distributions of $\chi_\phi^3 \nabla p \cdot \nabla p$ for the shielded structure (a) and the penetrated structure (b), respectively, where the diffusion coefficients are $\kappa^{\,\mathrm{p}}_\mathrm{solid} = 1$ and $\kappa^{\,\mathrm{p}}_\mathrm{void} = 1000$.
It can be seen that the value of $\chi_\phi^3 \nabla p \cdot \nabla p$ is larger for the shielded structure compared to the penetrated structure. Especially, the value is larger in the narrower solid domain.
Therefore, the penetrating effect on the design objects is enhanced by incorporating the evaluation function $J_p$ into the objective function of the minimization problem.

The sensitivity $J_p'$ of the evaluation function (\ref{eq:Jp}) is formulated as follows:
\begin{equation} \label{eq:Jp'}
J_p'=- 3 \chi_\phi^2 \nabla p \cdot \nabla p - \frac{\partial a_p}{\partial \chi_\phi} \nabla p \cdot \nabla {\hat{p} }
\end{equation}
where the adjoint variable $\hat{p} \in H^{1}(D)$ is defined as follows:
\begin{align} \label{eq:AdjP}
    \begin{cases}
         -\nabla \cdot (a_p \nabla \hat{p}) = 2 \chi_\phi^3 \nabla ^2 p  \, \, \, \, & \text{in} \, \, \, \, D\\
        \hat{p} = 0 \, \, \, \, & \text{on} \, \, \, \, \Gamma^{\,\mathrm{p}}_\mathrm{in}\\
        -\bm{n} \cdot \left(a_p \nabla \hat{p} \right) = 0 \qquad & \text{on} \, \, \, \,  \partial D \setminus \Gamma^{\,\mathrm{p}}_\mathrm{in}
    \end{cases}
\end{align}
\color{black}
\section{Application to Optimization Problem}
\subsection{Formulation of Optimization Problem Considering Target Geometric Features} \label{sec:formulation}
\color{black}
In this section, the minimum mean compliance problem with geometric and volume conditions is addressed as a concrete optimization problem and formulated as topology optimization based on the level set method.
We consider the domain $\Omega$ filled with linear isotropic homogeneous solid within the fixed design domain $D$, where the displacement $\bm{u}$ are fixed on the boundary $\Gamma_u$ and traction $\bm{t}$ is applied on the boundary $\Gamma_t$.
The upper limit of the volume is $V_{\max}$, and shielding and penetrating conditions are also included.
\red{The shielding and penetrating conditions are added as terms in the objective function rather than formulated as constraint functions. This is because their geometric effects are intended to be controlled by weighting parameters.}
Under these conditions, the optimization problem is formulated as follows:
\begin{alignat}{2}
  &\inf_{\chi_\phi }&
	&J_{u,s,p}=\: \displaystyle (1-\omega_s-\omega_p) J_u + \omega_s J_s + \omega_p J_p\\
  &\text{where}& &J_{u}=\: \displaystyle \int_{\Gamma_t} \bm{t}\cdot \bm{u} \dG\\
  &&&J_{s}=\int_D \left(1-\chi_{\phi}\right)^3 \, \nabla s \cdot \nabla s \dO \notag\\
  &&&J_{p}=\int_D \chi_{\phi}^3 \nabla p \cdot \nabla p \dO \notag\\
  &\text{subject to:}&\qquad &\int_{\Gamma_t} \bm{t} \cdot \tilde{\bm{u}} \dG
	-\int_D \bm{\epsilon}(\bm{u}): \chi_\phi\DD :\bm{\epsilon}(\tilde{\bm{u}}) \dO=0\\
  &&&\qquad \text{for  }\, \, \forall \tilde{\bm{u}} \in \U, \bm{u} \in \U \notag\\
  &&&\int_D a_s \nabla s \cdot \nabla {\tilde{s} } \dO =0\\
  &&&\qquad \text{for  }\, \,  \forall \tilde{s} \in \S, s \in \S \notag\\
  &&&\int_D a_p \nabla p \cdot \nabla {\tilde{p} } \dO =0\\
  &&&\qquad \text{for  }\, \,  \forall \tilde{p} \in \P, p \in \P \notag\\
  &&&\int_D \chi_{\phi} \dO -V_{\max} \le 0
\end{alignat}
where $\omega_s$ and $\omega_p$ are positive weighting parameters for the shielding and penetrating evaluation functions, and satisfy the relation $\omega_s+\omega_p<1$. Additionally, $\DD$ and $\bm{\epsilon}$ denote the elastic tensor and linearized strain tensor, respectively. The functional spaces 
$\U$, $\S$ and $\P$ are defined as follows.
\begin{align}
	&	\U= \left\{ \tilde{\bm{u}} \in H^1(D)^d \,\,  \Bigl |  \,\,\:\:\: \qquad \qquad \tilde{\bm{u}}  =0\quad {\rm on } \,\,\; \Gamma_u    \right\}\\
	&	\S= \left\{ \tilde{s} \in H^1(D) \quad \Bigl |  \,\,\:\:\:\, \qquad \qquad \tilde{s}  =1\quad {\rm on } \,\,\; \Gamma^{\,\mathrm{s}}_\mathrm{out} \,, \quad \tilde{s}  =0\quad {\rm on }  \,\,\; \Gamma^{\,\mathrm{s}}_\mathrm{in}    \right\}\\
        &	\P= \left\{ \tilde{p} \in H^1(D) \:\:\:\; \Bigl |   \,\, -\bm{n} \cdot \left(a_p\nabla \tilde{p}\right) =1 \quad {\rm on }\,\,\; \Gamma^{\,\mathrm{p}}_\mathrm{out} \,, \quad \tilde{p}  =0\quad {\rm on }  \,\,\; \Gamma^{\,\mathrm{p}}_\mathrm{in}    \right\}
\end{align}

\red{Regarding the penetrating condition, one possible approach is to directly define through holes as non-design domains before optimization.  
However, the proposed method also optimizes the through holes.
This optimization has the great advantage of eliminating the effort required to define the penetrating region and minimizing the deterioration of mean compliance.}

The sensitivity of the objective function $J_{u,s,p}$ is
\begin{equation}
J_{u,s,p}' =\: \displaystyle (1-\omega_s-\omega_p) J_u' + \omega_s J_s' + \omega_p J_p'
\end{equation}
$J_s'$ and $J_p'$ are as given by equations (\ref{eq:Js'}) and (\ref{eq:Jp'}). And the sensitivity of the minimum mean compliance problem in three dimensions  \cite{bonnet2013topological,otomori2015matlab} is given by the following equation as the topological derivative.
\begin{equation}
J_u'=\chi_\phi \bm{\epsilon}(\bm{u}): \mathbb{A}  :\bm{\epsilon}(\tilde{\bm{u}})
\end{equation}
The coefficient tensor $\mathbb{A}$ is 
\begin{equation} \label{eq:tensor A}
\mathbb{A}_{ijkl}=\frac{3(1-\nu)}{2(1+\nu)(7-5\nu)}\left[
-\frac{(1-14\nu+15\nu^2)E}{(1-2\nu)^2}\delta_{ij}\delta_{kl}
+5E(\delta_{ik}\delta_{jl}+\delta_{il}\delta_{jk})
\right]
\end{equation}
where $E$, $\nu$ and $\delta_{ij}$ denote the Young's modulus, Poisson ratio and Kronecker delta, respectively.
Equation (\ref{eq:tensor A}) is for the three-dimensional case, but it can also be applied to two-dimensional problems by assuming plane stress conditions.
Since the minimum mean compliance problem is a self-adjoint problem, the adjoint variable with respect to the displacement field is identical to the displacement field itself.

In this section, the formulation considers both shielding and penetrating conditions. However, please note that if either the shielding condition or the penetrating condition is not required, the corresponding equations need not be considered.
\color{black}
\subsection{Optimization Algorithm} \label{sec:flowchart}
\begin{figure}[H]
	\begin{center}
		\includegraphics[height=7cm]{ 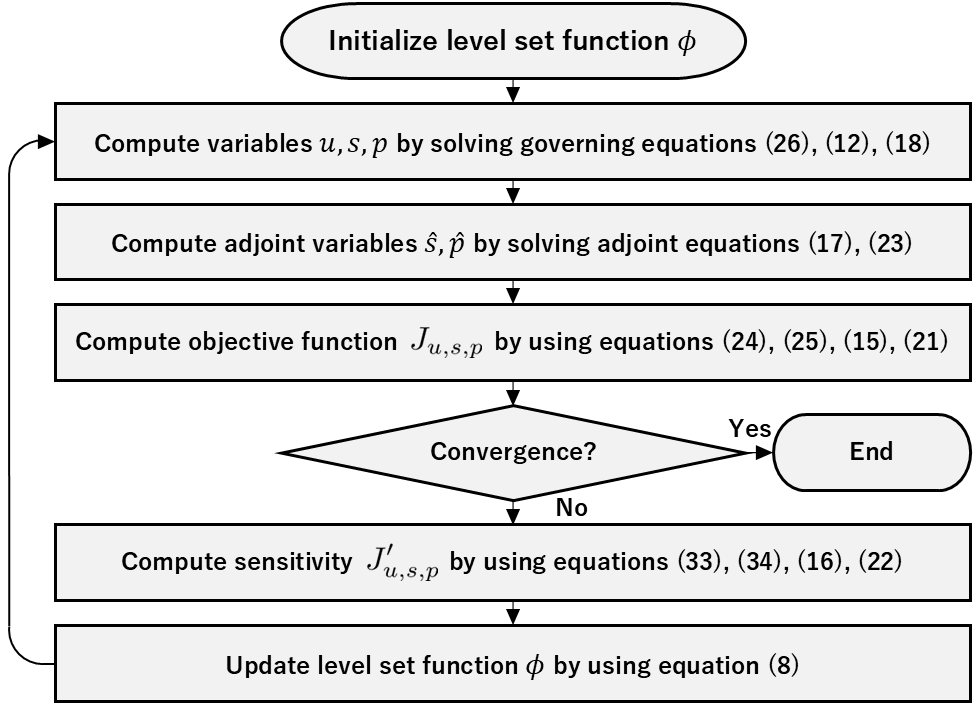}
		\caption{Flowchart of the optimization algorithm.}
		\label{fig:flowchart}
	\end{center}
\end{figure}

Here, we describe the optimization algorithm.
The volume is constrained by the augmented Lagrangian method.
Figure \ref{fig:flowchart} shows the optimization algorithm of the proposed method.
In this paper, the initial level set function $\phi$ is set to 1 for the whole domain, and solving equations are used the finite element method.
In the presented optimization examples, the domain is discretized into triangular elements in 2D and tetrahedral elements in 3D, but since the fictitious physical equations are linear isotropic diffusion equations, specialized numerical techniques for mesh discretization are not required.

\color{black}
\section{Numerical Examples}
In this section, we validate the proposed methods through  numerical examples.
First, we assess the proposed method for shielding through examples~\ref{sec:paramS effect},~\ref{sec:2D-S}, and~\ref{sec:3D-S}. Next, we analyze the proposed method for penetrating through examples~\ref{sec:paramP effect},~\ref{sec:2D-P}, and~\ref{sec:3D-P}.  Finally, we present the optimization example that considers both shielding and penetrating conditions in \red{examples~\ref{sec:2D-SP},~\ref{sec:3D-SP}}.
All of these numerical examples are derived from the minimum mean compliance problem formulated in Section~\ref{sec:formulation}, with the Young's modulus $E$ set to $210 \times 10^9$ and the Poisson ratio $\nu$ set to $0.3$.
Additionally, traction $\bm{t}$ is applied downward in all cases, with the magnitude of $1.0 \times 10^5$ in 2D and $1.0 \times 10^2$ in 3D.
Furthermore, for all examples, the transition width $\delta$ in the equation (\ref{eq:chi reg}) is set to $1.0 \times 10^{-3}$, and the regularization parameter $\tau$ in the equation (\ref{eq:reaction_diffusion_eq}) is generally set to $1.0 \times 10^{-4}$ unless stated otherwise.
It should be noted that the finite element analysis for all computational examples is performed using FreeFEM++ software\cite{hecht2012new}, and the visualization of the data is carried out using the open-source software ParaView\cite{ahrens2005paraview}.
\subsection{Effect of Parameters on Shielding Feature} \label{sec:paramS effect}
\begin{figure}[H]
    \begin{center}
        \begin{minipage}[t]{0.49\columnwidth}
            \centering
            \includegraphics[width=0.95\columnwidth]{ 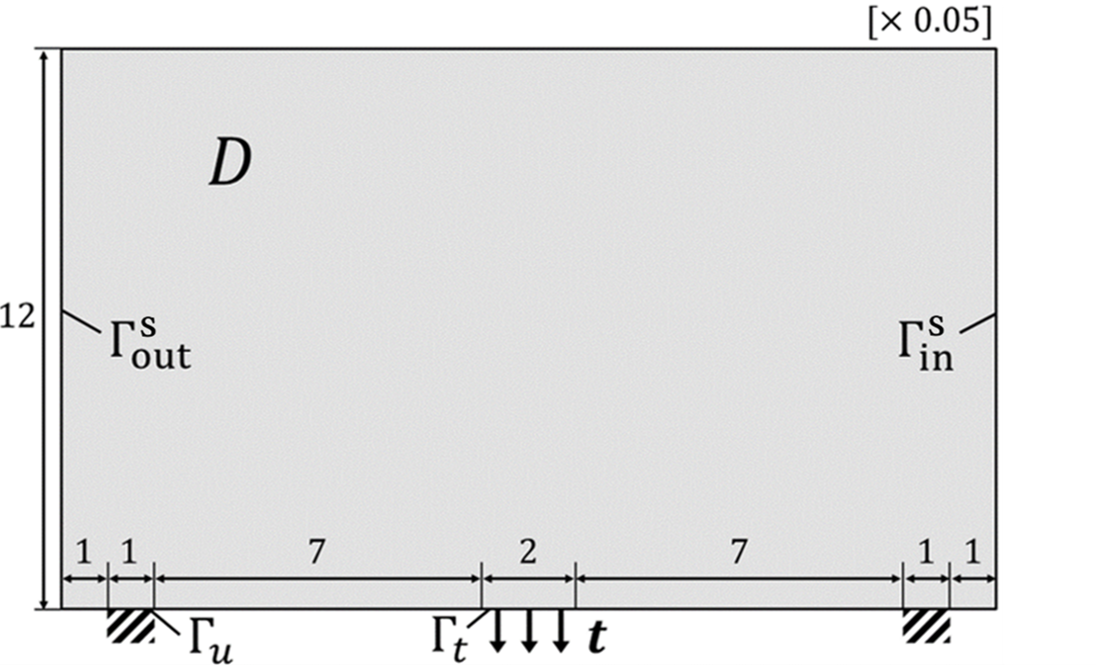}
            \subcaption{Problem setting for mean compliance minimization considering the shielding condition}
        \end{minipage}
        \begin{minipage}[t]{0.49\columnwidth}
            \centering
            \includegraphics[width=0.95\columnwidth]{ 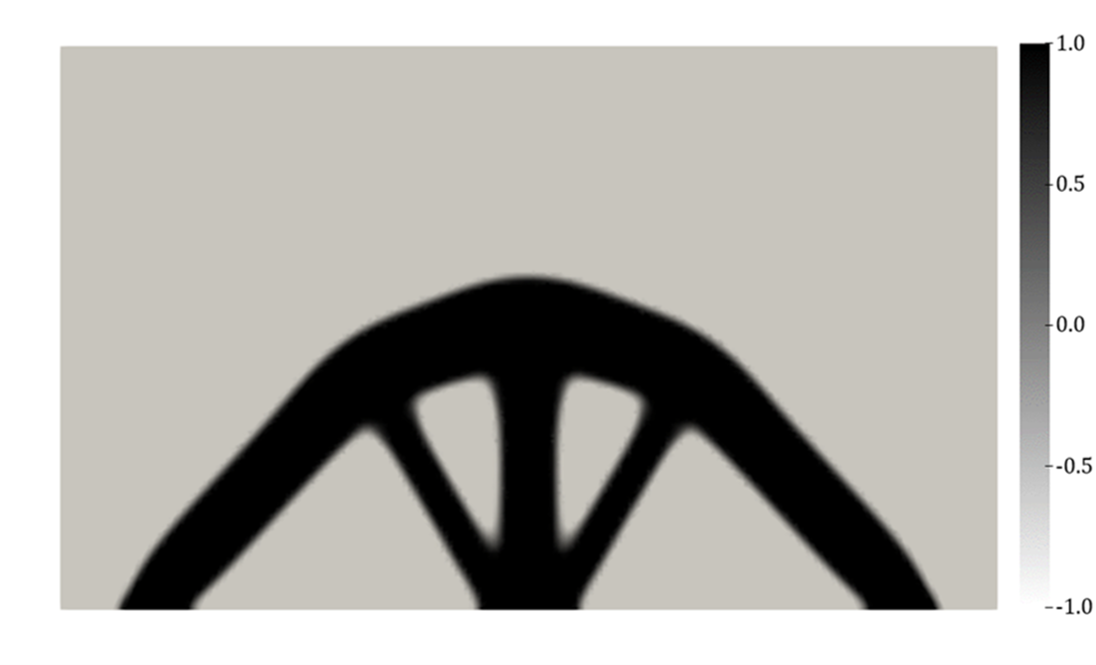}
            \subcaption{Configuration of optimized $\phi$ without the shielding condition}
        \end{minipage}
        \caption{Target for comparing the effect of parameters under the shielding condition.}
	\label{fig:paramS}
    \end{center}
\end{figure}

First, we consider the effect of the parameters in the proposed method for shielding.
The parameters of interest are the ratio of the diffusion coefficients $\kappa^{\,\mathrm{s}}_\mathrm{solid}$ and $\kappa^{\,\mathrm{s}}_\mathrm{void}$, as well as the weighting parameter $\omega_s$.
As the setting for this analysis, $\kappa^{\,\mathrm{s}}_\mathrm{solid}$ is fixed at $1$.
The reason for this setting is that the fictitious field $s$ depends on the ratio of the diffusion coefficients, not their magnitudes, since the governing equation (\ref{eq:GovS}) does not include source terms.
The problem setup for examining the influence of parameters is illustrated in Figure \ref{fig:paramS}(a). The displacement is fixed at the left and right boundaries $\Gamma_u$ on the lower edge, and the traction $\bm{t}$ is applied to the center boundary $\Gamma_t$ on the lower edge. As boundary conditions for the fictitious physical field $s$, the left edge and the right edge are defined as the boundaries $\Gamma^{\,\mathrm{s}}_\mathrm{out}$ and $\Gamma^{\,\mathrm{s}}_\mathrm{in}$, respectively. Additionally, the maximum volume is limited to 20\% of the fixed design domain $D$.
The optimized result without considering the shielding condition is shown in Figure \ref{fig:paramS}(b), and the optimized structure does not have any shielding walls.

\begin{figure}[H]
	\begin{center}
		\includegraphics[width=12cm]{ 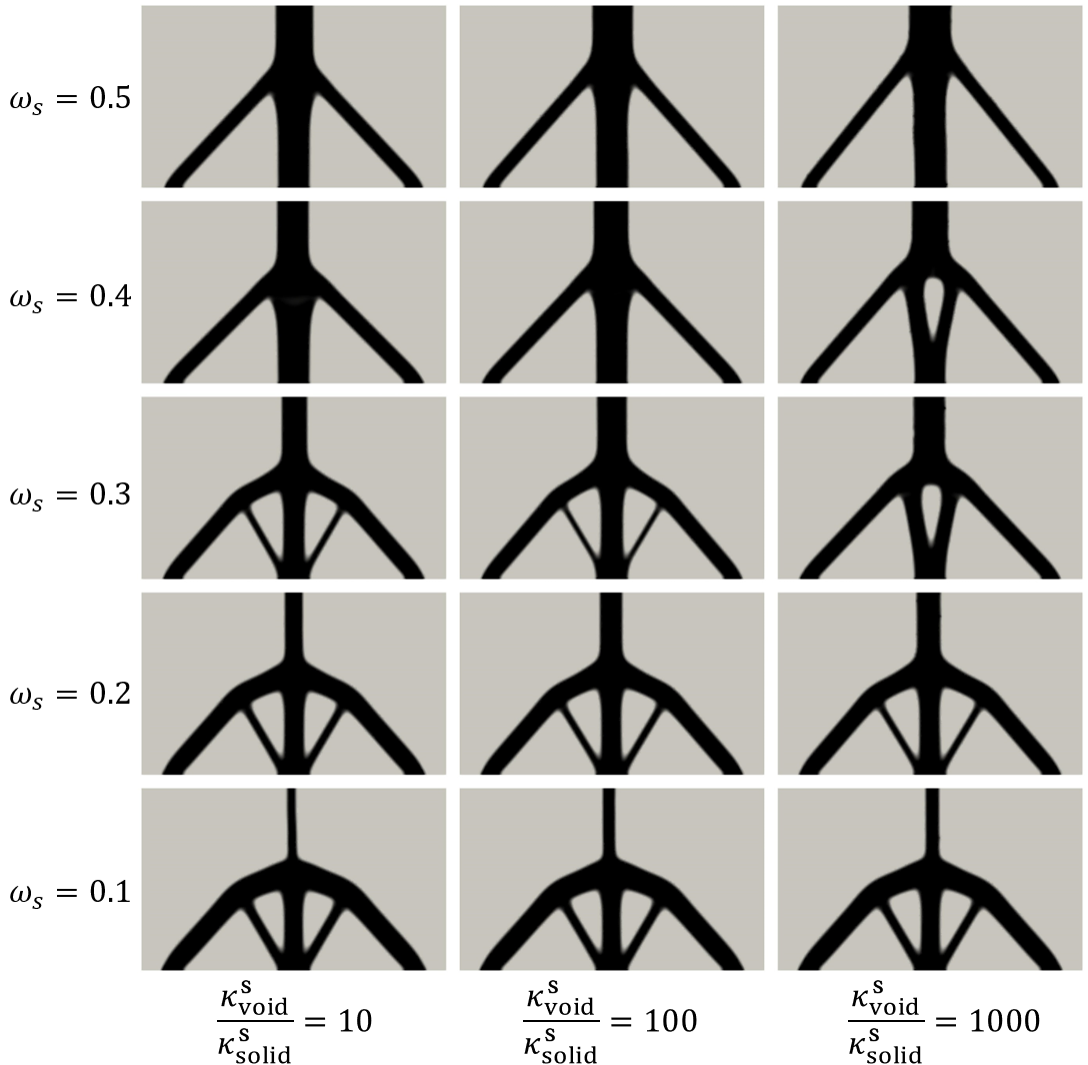}
		\caption{Effect of diffusion coefficient ratio $\displaystyle \frac{\kappa^{\,\mathrm{s}}_\mathrm{void}}{\kappa^{\,\mathrm{s}}_\mathrm{solid}}$ and weighting parameter $\omega_s$ on optimization considering the shielding condition.}
		\label{fig:tableS}
	\end{center}
\end{figure}

We investigate how parameters for shielding affect geometric features of the target optimization.
Specifically, we examine 15 patterns, where $\kappa^{\,\mathrm{s}}_\mathrm{solid}$ takes three values of (10, 100, 1000) for the ratio of the diffusion coefficients, and $\omega_s$ takes five values of (0.1, 0.2, 0.3, 0.4, 0.5).
The results are as depicted in Figure \ref{fig:tableS}.
It is confirmed that all structures include shielding walls.
Regarding the relationship between the diffusion coefficient ratio and the thickness of the shielding walls, it is observed that the thickness is larger for $\kappa^{\,\mathrm{s}}_\mathrm{void}=100$ compared to $\kappa^{\,\mathrm{s}}_\mathrm{void}=10$. However, the thickness is almost the same for $\kappa^{\,\mathrm{s}}_\mathrm{void}=100$ and $\kappa^{\,\mathrm{s}}_\mathrm{void}=1000$, while there are some differences in shape.
From this result, it is desirable to set $\kappa^{\,\mathrm{s}}_\mathrm{void}$ to approximately $100$ or higher for effective shielding. 
Therefore, in the following numerical examples, $\kappa^{\,\mathrm{s}}_\mathrm{solid}$ and $\kappa^{\,\mathrm{s}}_\mathrm{void}$ are set to $1$ and $100$, respectively.
Next, regarding the relationship between the weighting parameter and the thickness of the shielding walls, it is evident that as the value of $\omega_s$ increases, the thickness becomes noticeably larger.
This result indicates that the shielding effectiveness can be controlled by varying the weighting parameter $\omega_s$.
\subsection{2D Example Under Shielding Condition} \label{sec:2D-S}
\begin{figure}[H]
	\begin{center}
		\includegraphics[height=5cm]{ 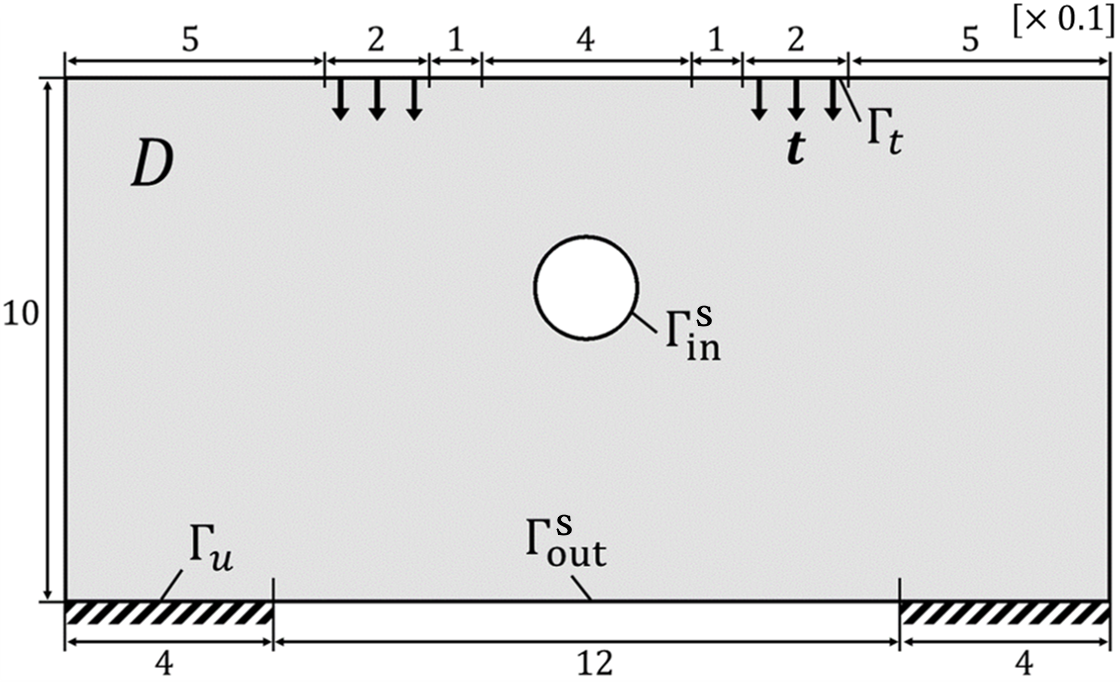}
		\caption{Problem setting for 2D mean compliance minimization considering the shielding condition.}
		\label{fig:2Sset}
	\end{center}
\end{figure}

As the second validation of the proposed method for shielding, we examine the two-dimensional optimization examples and discuss their fictitious physical fields.
The problem setup is shown in Figure \ref{fig:2Sset}. The displacement is fixed at the left and right boundaries $\Gamma_u$ on the lower edge, and the traction $\bm{t}$ is applied from the two upper boundaries $\Gamma_t$.
In addition, to control shielding, the boundary $\Gamma^{\,\mathrm{s}}_\mathrm{out}$ is defined at the center of the lower edge, and the boundary $\Gamma^{\,\mathrm{s}}_\mathrm{in}$ is set on the perimeter of the central circular non-design domain.
Furthermore, the maximum volume is limited to $50\%$ of the fixed design domain $D$, and the weighting parameter $\omega_s$ is set to $0.06$.

\begin{figure}[H]
    \centering
    \begin{minipage}[t]{0.46\linewidth}
        \centering
        \includegraphics[width=\linewidth]{ 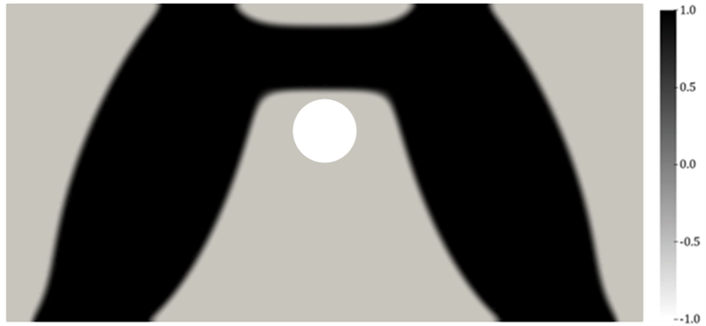}
        \subcaption{Configuration of optimized $\phi$ without the shielding condition}
    \end{minipage}
    \hspace{0.03\linewidth}
    \begin{minipage}[t]{0.46\linewidth}
        \centering
        \includegraphics[width=\linewidth]{ 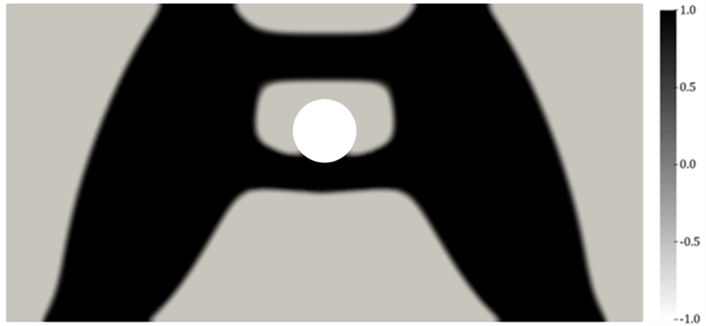}
        \subcaption{Configuration of optimized $\phi$ with the shielding condition}
    \end{minipage}
    \vspace{0.02\linewidth}
    \begin{minipage}[t]{0.46\linewidth}
        \centering
        \includegraphics[width=\linewidth]{ 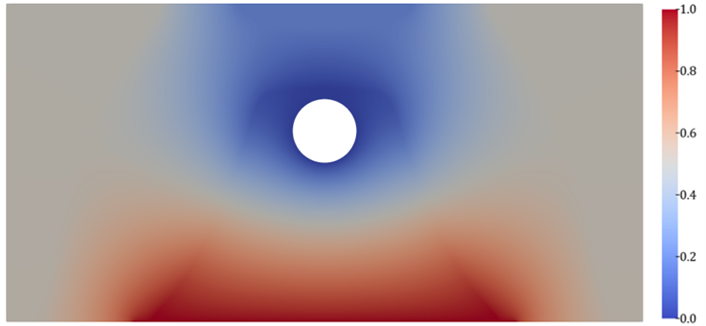}
        \subcaption{Fictitious physical field $s$ without the shielding condition}
    \end{minipage}
    \hspace{0.02\linewidth}
    \begin{minipage}[t]{0.46\linewidth}
        \centering
        \includegraphics[width=\linewidth]{ 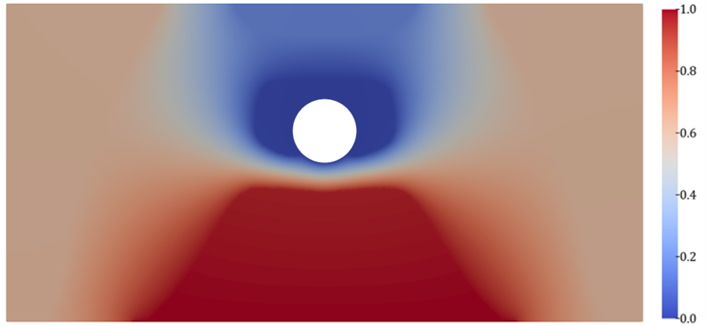}
        \subcaption{Fictitious physical field $s$ with the shielding condition.}
    \end{minipage}
    \caption{Comparison of optimized 2D results without and with the shielding condition.}
    \label{fig:2Sopt}
\end{figure}

(a) and (b) in Figure \ref{fig:2Sopt} represent optimized configuration under conditions without and with shielding, respectively.
It is clear that by introducing the shielding condition, the optimized structure shields between the boundaries $\Gamma^{\,\mathrm{s}}_\mathrm{out}$ and $\Gamma^{\,\mathrm{s}}_\mathrm{in}$.
The fictitious physical fields $s$ in these examples are as shown in (c) and (d) of Figure \ref{fig:2Sopt}.
Without the shielding condition, the temperature changes significantly in the void domain between the boundaries $\Gamma^{\,\mathrm{s}}_\mathrm{out}$ and $\Gamma^{\,\mathrm{s}}_\mathrm{in}$. In contrast, with the shielding condition, the shielded part is formed and the temperature change is concentrated in that area. As a result, the temperature change is nearly zero in the void domain.
\subsection{3D Example Under Shielding Condition} \label{sec:3D-S}
\begin{figure}[H]
	\begin{center}
		\includegraphics[height=5cm]{ 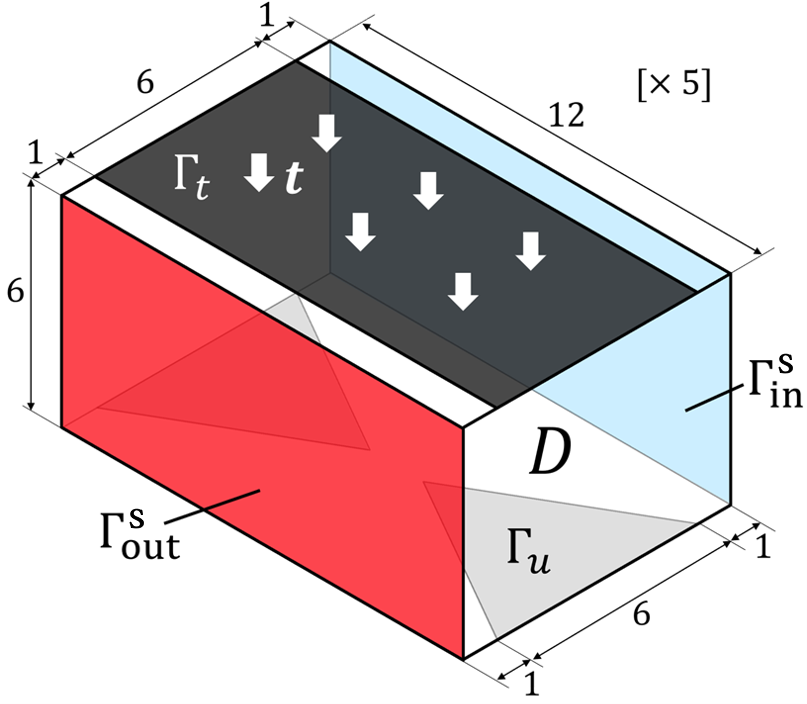}
		\caption{Problem setting for 3D mean compliance minimization considering the shielding condition.}
		\label{fig:3Sset}
	\end{center}
\end{figure}

In the third validation of the shielding effectiveness, we present three-dimensional optimization examples and review their convergence histories.
The problem setting is illustrated in Figure \ref{fig:3Sset}. 
In the fixed design domain $D$, the displacement is fixed at the boundaries $\Gamma_u$ of the symmetric equilateral triangular regions on the bottom surface, and the traction $\bm{t}$ is applied to a part of the top surface. Additionally, to obtain the shielded structure, the front and back surfaces are designated as the boundaries $\Gamma^{\,\mathrm{s}}_\mathrm{out}$ and $\Gamma^{\,\mathrm{s}}_\mathrm{in}$, respectively. Furthermore, the maximum volume is limited to $40\%$ of the fixed design domain $D$, and the weighting parameter $\omega_s$ is set to $0.1$.

\begin{figure}[H]
    \centering
    \begin{minipage}[t]{0.46\linewidth}
        \centering
        \includegraphics[width=\linewidth]{ 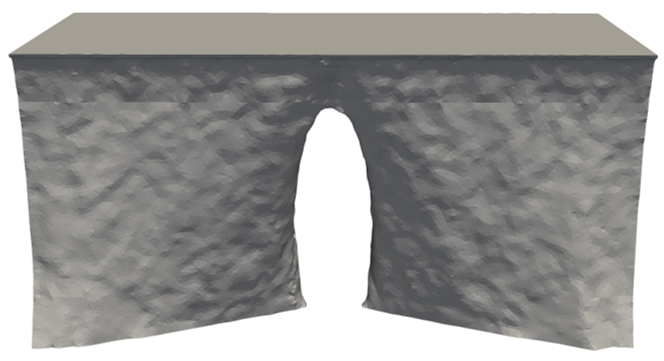}
        \subcaption{Front view without the shielding condition}
    \end{minipage}
    \hspace{0.03\linewidth}
    \begin{minipage}[t]{0.46\linewidth}
        \centering
        \includegraphics[width=\linewidth]{ 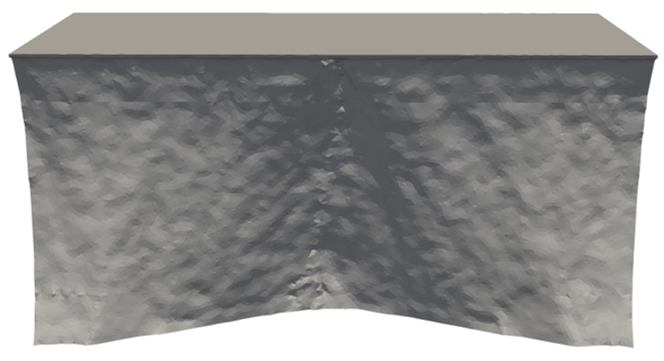}
        \subcaption{Front view with the shielding condition}
    \end{minipage} 
    \vspace{0.02\linewidth}
    \begin{minipage}[t]{0.46\linewidth}
        \centering
        \includegraphics[width=\linewidth]{ 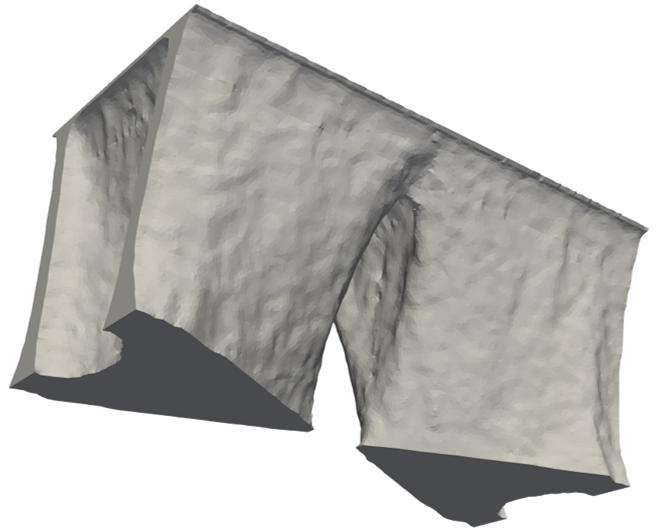}
        \subcaption{Oblique view without the shielding condition}
    \end{minipage}
    \hspace{0.02\linewidth}
    \begin{minipage}[t]{0.46\linewidth}
        \centering
        \includegraphics[width=\linewidth]{ 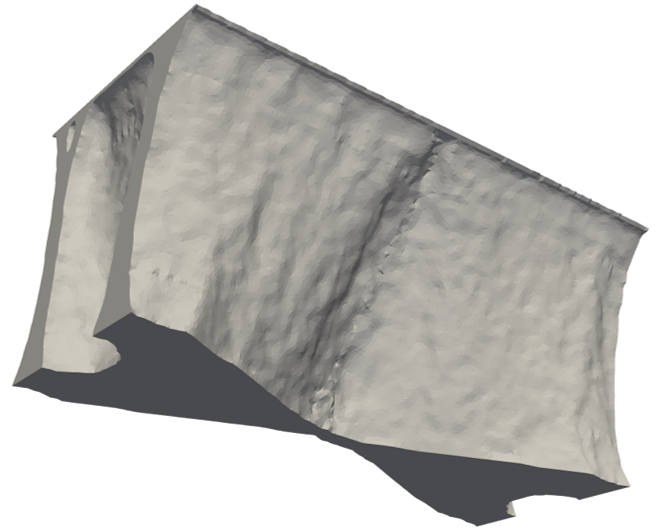}
        \subcaption{Oblique view with the shielding condition}
    \end{minipage}
    \caption{Comparison of optimized 3D configurations without and with the shielding condition.}
    \label{fig:3Sopt}
\end{figure}

Figure \ref{fig:3Sopt} compares the optimized structures without and with the shielding condition.
As shown in (a) and (c), in the case without the shielding condition, the optimized structure is penetrated as an arch.
However, as shown in (b) and (d), in the case with the shielding condition, the optimized structure forms a shielding wall that obstructs the area between the boundaries $\Gamma^{\,\mathrm{s}}_\mathrm{out}$ and $\Gamma^{\,\mathrm{s}}_\mathrm{in}$.
This result demonstrates that the proposed method is useful for shielding even in three-dimensional problems.

\begin{figure}[H]
    \begin{center}
        \begin{minipage}[t]{0.32\columnwidth}
            \centering
            \includegraphics[width=0.95\columnwidth]{ 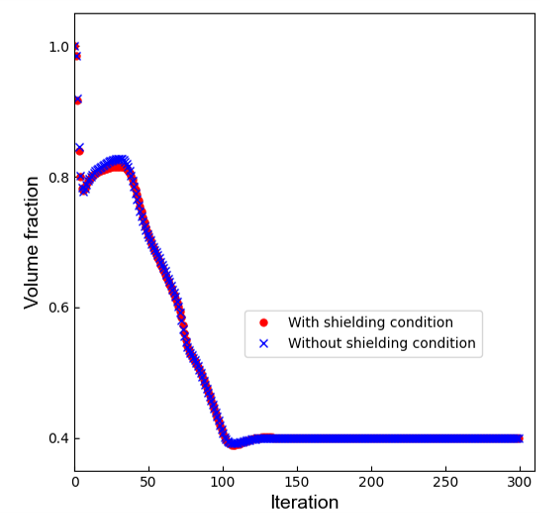}
            \subcaption{Change in volume fraction}
        \end{minipage}
        \begin{minipage}[t]{0.32\columnwidth}
            \centering
            \includegraphics[width=0.95\columnwidth]{ 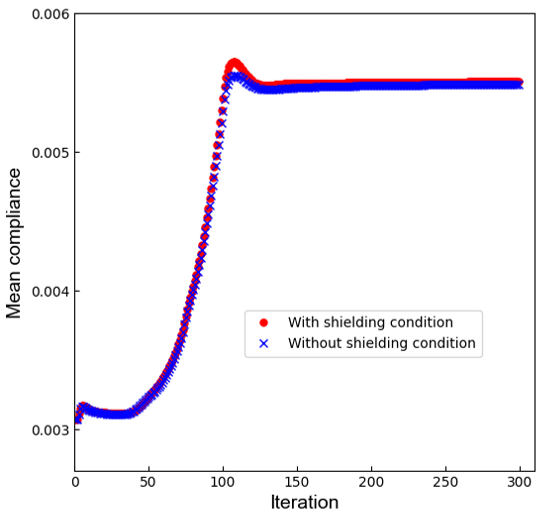}
            \subcaption{Change in mean compliance $J_u$}
        \end{minipage}
         \begin{minipage}[t]{0.32\columnwidth}
            \centering
            \includegraphics[width=0.95\columnwidth]{ 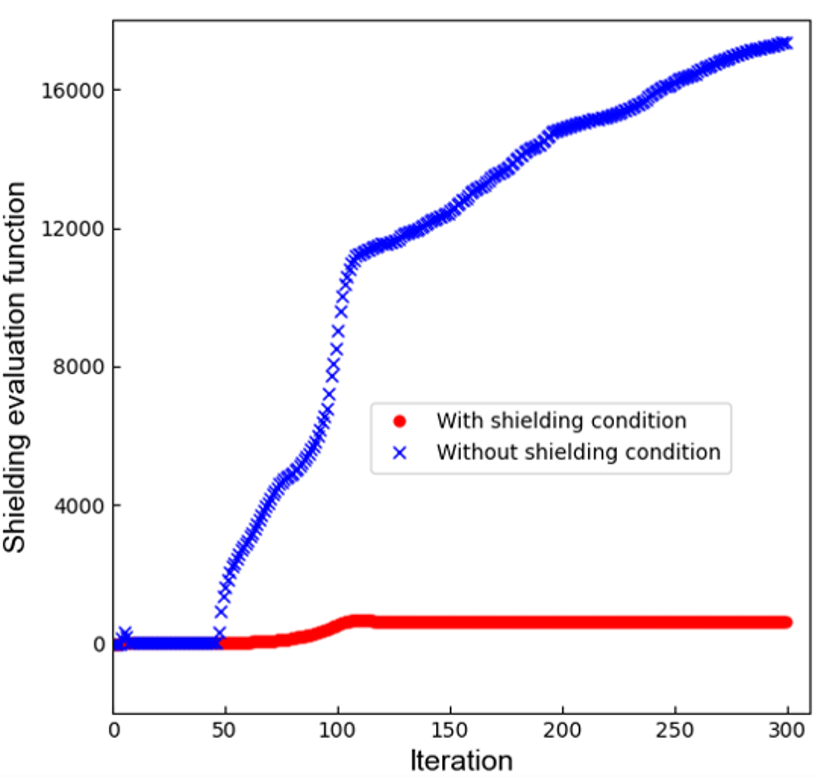}
            \subcaption{Change in shielding evaluation function $J_s$}
        \end{minipage}
        \caption{Convergence histories comparing optimization without and with the shielding condition.}
    \label{fig:graph3S}
    \end{center}
\end{figure}

Next, Figure \ref{fig:graph3S} shows the convergence histories of these optimization examples.
As shown in (a), both optimizations are implemented to converge the volume within approximately 100 iterations out of the total 300 iterations.
(b) illustrates the change in mean compliance for each case, and there is little difference between them.
In contrast, (c) shows the change in the values of the shielding evaluation function, and there is significant difference between the cases.
Specifically, without the shielding condition, the value increases significantly as the structure is penetrated. However, with the shielding condition, it works to decrease the value. As a result, the structure remains shielded.
\subsection{Effect of Parameters on Penetrating Feature} \label{sec:paramP effect}
As the first consideration for penetrating, we examine the effect of the parameters in the proposed method for penetrating. 
The parameters of interest are the ratio of the diffusion coefficients $\kappa^{\,\mathrm{p}}_\mathrm{solid}$ and $\kappa^{\,\mathrm{p}}_\mathrm{void}$, as well as the weighting parameter $\omega_p$.
As the setting for this analysis, $\kappa^{\,\mathrm{p}}_\mathrm{solid}$ is fixed at $1$.
The reason for this setting is that the shape of the fictitious field $p$ depends on the ratio of the diffusion coefficients, not their magnitudes.
In the governing equation (\ref{eq:GovP}), the Neumann boundary condition with a constant value of 1 is applied at the boundary $\Gamma^{\,\mathrm{p}}_\mathrm{out}$. As a result, the scale of the fictitious physical field changes depending on the magnitude of the diffusion coefficients.
However, as shown in the equation (\ref{eq:reaction-diffusion}), since normalized sensitivity is used to update the design variable, the scale of the fictitious physical field does not affect the optimization results for the model used in this paper. Therefore, the magnitude of the diffusion coefficients does not influence the effectiveness for penetrating.

\begin{figure}[H]
    \begin{center}
        \begin{minipage}[t]{0.49\columnwidth}
            \centering
            \includegraphics[width=0.95\columnwidth]{ 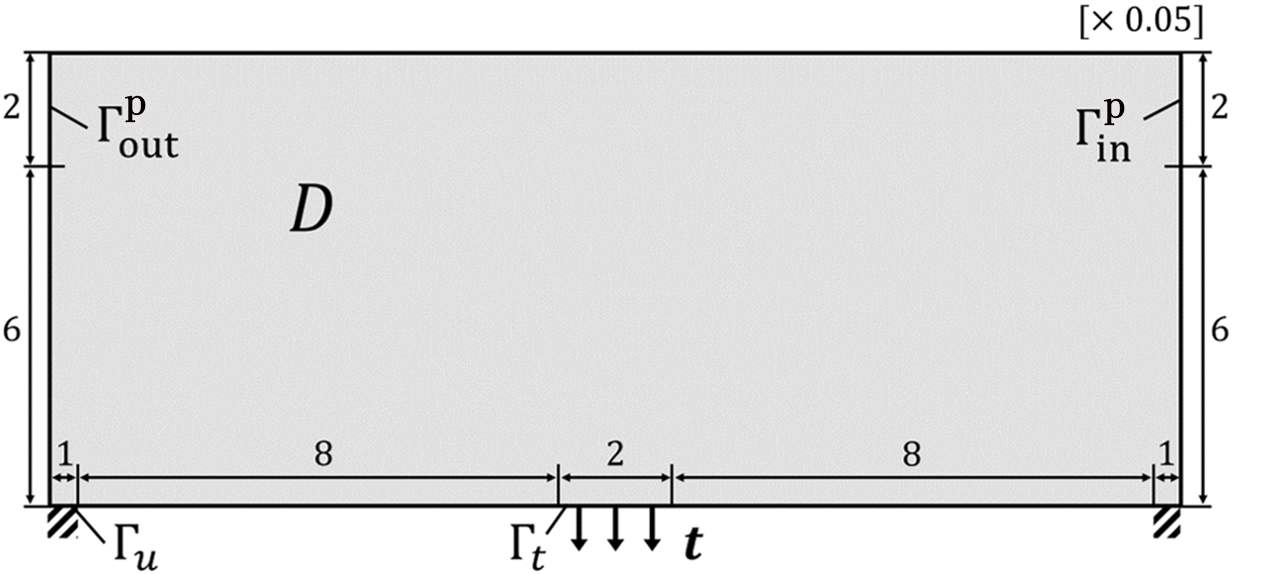}
            \subcaption{Problem setting for mean compliance minimization considering the penetrating condition}
        \end{minipage}
        \begin{minipage}[t]{0.49\columnwidth}
            \centering
            \includegraphics[width=0.95\columnwidth]{ 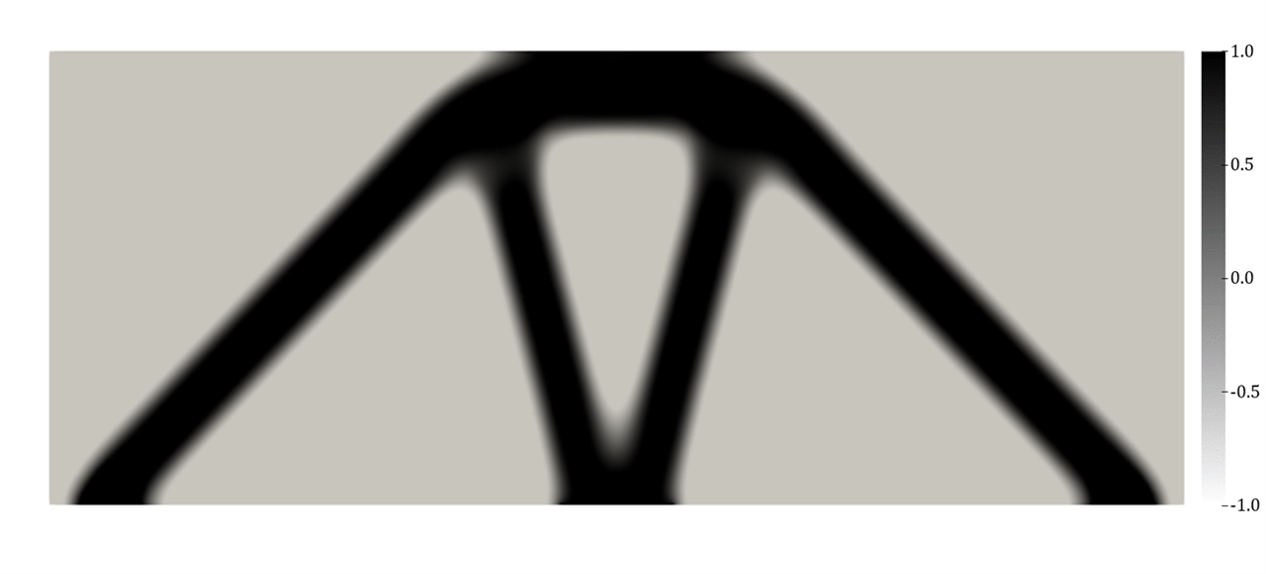}
            \subcaption{Configuration of optimized $\phi$ without the shielding condition}
        \end{minipage}
        \caption{Target for comparing the effect of parameters under the penetrating condition.}
	\label{fig:paramP}
    \end{center}
\end{figure}

The problem setup for examining the influence of parameters is illustrated in Figure \ref{fig:paramP}(a). The displacement is fixed at the left and right boundaries $\Gamma_u$ on the lower edge, and the traction $\bm{t}$ is applied to the center boundary $\Gamma_t$ on the lower edge. As boundary conditions for the fictitious physical field $p$, The upper parts of the left and right edges are defined as the boundaries \( \Gamma^{\,\mathrm{p}}_\mathrm{out} \) and \( \Gamma^{\,\mathrm{p}}_\mathrm{in} \), respectively.
Additionally, the maximum volume is limited to 25\% of the fixed design domain $D$.
The optimized result without considering the penetrating condition is shown in Figure \ref{fig:paramP}(b), and the optimized structure is not penetrated.

\begin{figure}[H]
	\begin{center}
		\includegraphics[width=12cm]{ 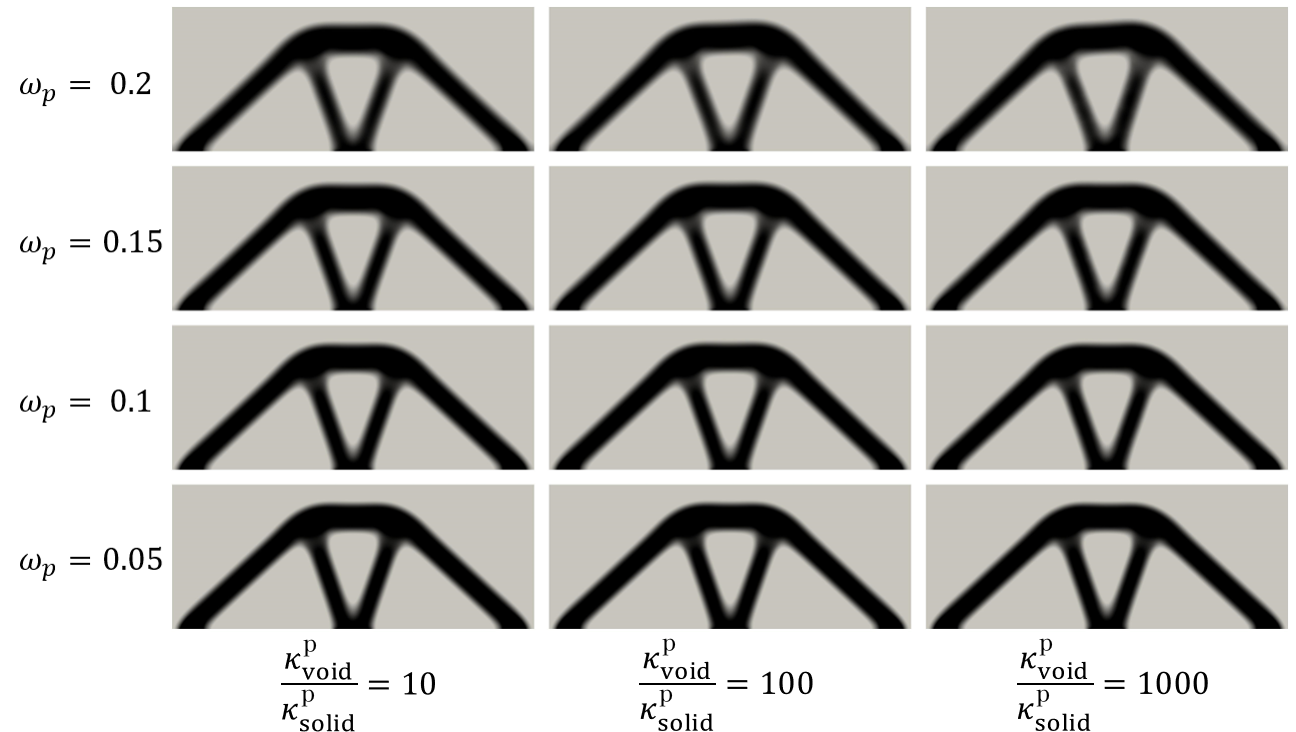}
		\caption{Effect of diffusion coefficient ratio $\displaystyle \frac{\kappa^{\,\mathrm{p}}_\mathrm{void}}{\kappa^{\,\mathrm{p}}_\mathrm{solid}}$ and weighting parameter $\omega_p$ on optimization considering the penetrating condition.}
		\label{fig:tableP}
	\end{center}
\end{figure}

We investigate how parameters for penetrating affect geometric features of the target optimization.
Specifically, we examine 12 patterns, where $\kappa^{\,\mathrm{p}}_\mathrm{solid}$ takes three values of (10, 100, 1000) for the ratio of the diffusion coefficients, and $\omega_p$ takes four values of (0.05, 0.1, 0.15, 0.2).
The results are as depicted in Figure \ref{fig:tableP}.
It is confirmed that all structures are penetrated at the top.
However, these structures show no noticeable geometric differences by changing the values of the diffusion coefficient ratio and the weighting parameter.
This result indicates that the penetrating effectiveness is difficult to control by varying the tested parameters, and it remains the subject for future research.
Based on this verification, in the following numerical examples, $\kappa^{\,\mathrm{p}}_\mathrm{solid}$ and $\kappa^{\,\mathrm{p}}_\mathrm{void}$ are set to $1$ and $100$, respectively.
\subsection{2D Example Under Penetrating Condition} \label{sec:2D-P}
\begin{figure}[H]
	\begin{center}
		\includegraphics[height=5cm]{ 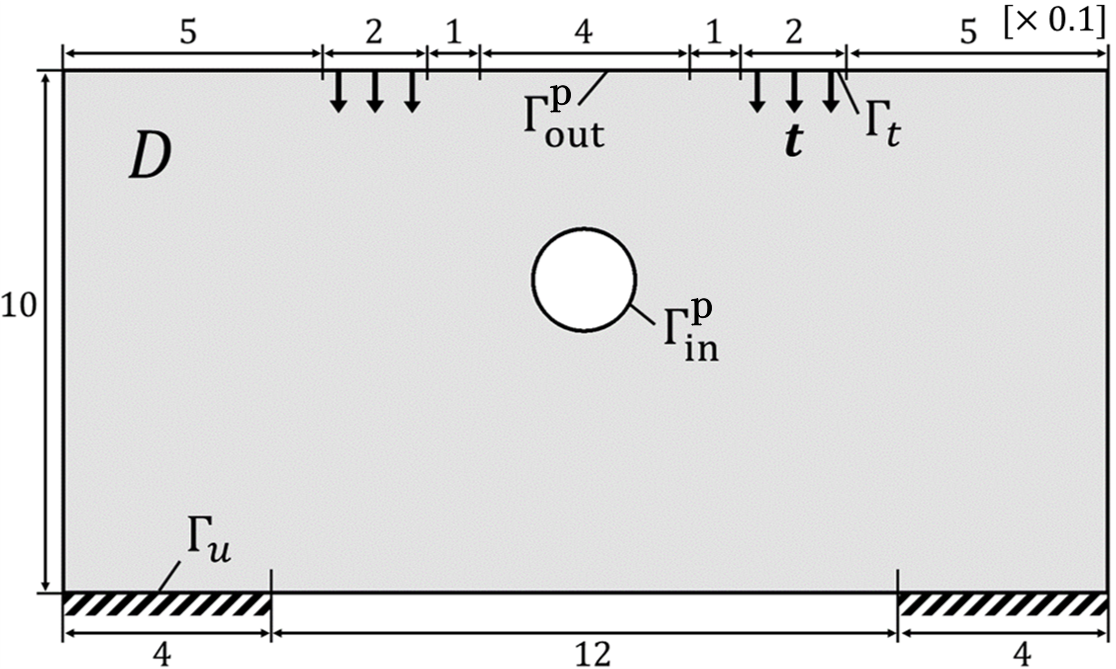}
		\caption{Problem setting for 2D mean compliance minimization considering the penetrating condition.}
		\label{fig:2Pset}
	\end{center}
\end{figure}

As the second validation of the proposed penetrating method, we examine the two-dimensional optimization examples and discuss their fictitious physical fields.
The problem setup is shown in Figure \ref{fig:2Pset}. Similar to Section~\ref{sec:2D-S}, the displacement is fixed at the left and right boundaries $\Gamma_u$ on the lower edge, and the traction $\bm{t}$ is applied from the two upper boundaries $\Gamma_t$.
In addition, to control penetrating, the boundary $\Gamma^{\,\mathrm{p}}_\mathrm{out}$ is defined at the center of the upper edge, and the boundary $\Gamma^{\,\mathrm{p}}_\mathrm{in}$ is set on the perimeter of the central circular non-design domain.
Furthermore, the maximum volume is limited to $50\%$ of the fixed design domain $D$, and the weighting parameter $\omega_p$ is set to $0.06$.

\begin{figure}[H]
    \centering
    \begin{minipage}[t]{0.46\linewidth}
        \centering
        \includegraphics[width=\linewidth]{ 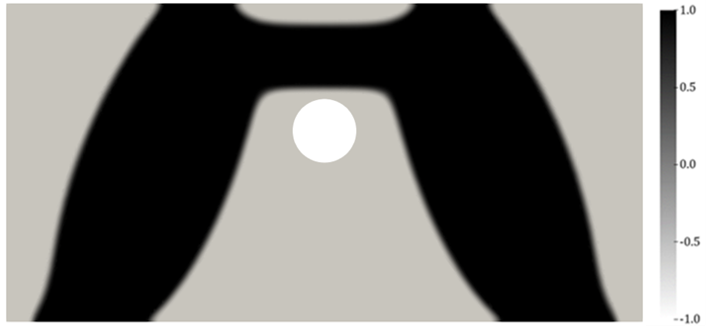}
        \subcaption{Configuration of optimized $\phi$ without the penetrating condition}
    \end{minipage}
    \hspace{0.03\linewidth}
    \begin{minipage}[t]{0.46\linewidth}
        \centering
        \includegraphics[width=\linewidth]{ 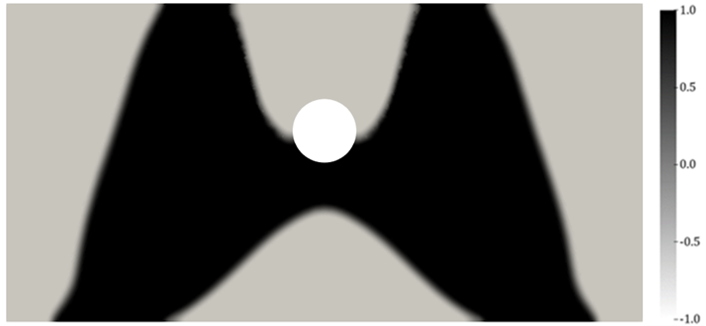}
        \subcaption{Configuration of optimized $\phi$ with the penetrating condition}
    \end{minipage}
    \vspace{0.02\linewidth}
    \begin{minipage}[t]{0.46\linewidth}
        \centering
        \includegraphics[width=\linewidth]{ 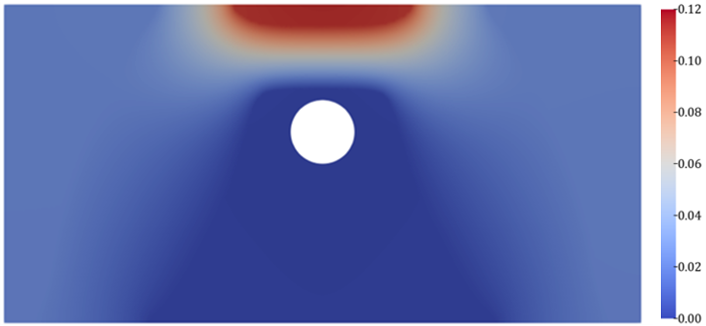}
        \subcaption{Fictitious physical field $p$ without the penetrating condition}
    \end{minipage}
    \hspace{0.02\linewidth}
    \begin{minipage}[t]{0.46\linewidth}
        \centering
        \includegraphics[width=\linewidth]{ 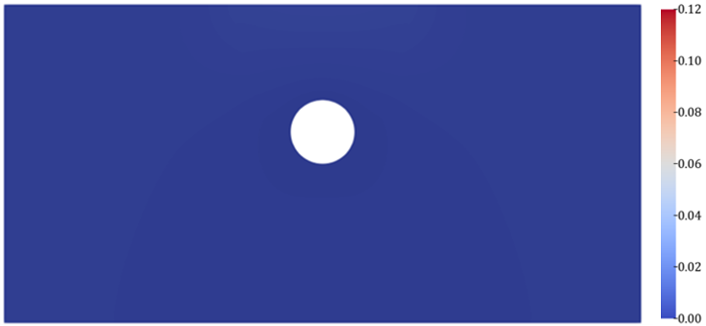}
        \subcaption{Fictitious physical field $p$ with the penetrating condition}
    \end{minipage}
    \caption{Comparison of optimized 2D results without and with the penetrating condition.}
    \label{fig:2Popt}
\end{figure}

(a) and (b) in Figure \ref{fig:2Popt} represent optimized configuration under conditions without and with penetrating, respectively.
It is clear that by introducing the penetrating condition, the optimized structure is penetrated between the boundaries $\Gamma^{\,\mathrm{p}}_\mathrm{out}$ and $\Gamma^{\,\mathrm{p}}_\mathrm{in}$.
The fictitious physical fields $p$ in these examples are as shown in (c) and (d) of Figure \ref{fig:2Popt}.
Without the penetrating condition, the temperature changes significantly in the solid domain between the boundaries $\Gamma^{\,\mathrm{p}}_\mathrm{out}$ and $\Gamma^{\,\mathrm{p}}_\mathrm{in}$. In contrast, with the penetrating condition, the shielded part is removed, resulting in the nearly uniform temperature overall. As a result, the temperature change is nearly zero in the solid domain.
\subsection{3D Example Under Penetrating Condition} \label{sec:3D-P}
\begin{figure}[H]
	\begin{center}
		\includegraphics[height=5cm]{ 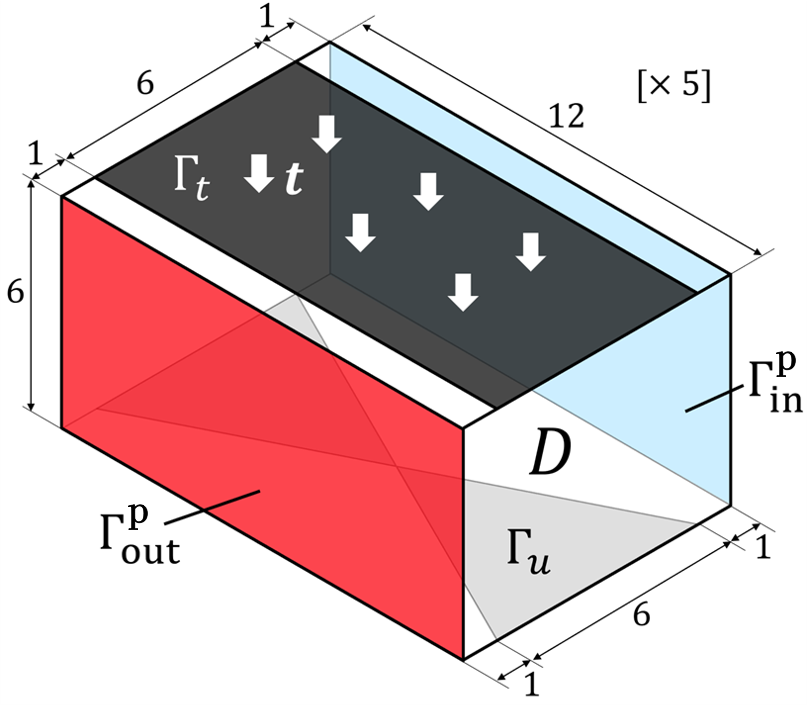}
		\caption{Problem setting for 3D mean compliance minimization considering the penetrating condition.}
		\label{fig:3Pset}
	\end{center}
\end{figure}

In the third validation of the penetrating effectiveness, we present three-dimensional optimization examples and review their convergence histories.
The problem setting is illustrated in Figure \ref{fig:3Pset}. 
In the fixed design domain $D$, the displacement is fixed at the boundaries $\Gamma_u$ of the adjacent symmetric isosceles triangular regions on the bottom surface, and the traction $\bm{t}$ is applied to a part of the top surface. Additionally, to obtain the penetrated structure, the front and back surfaces are designated as the boundaries $\Gamma^{\,\mathrm{p}}_\mathrm{out}$ and $\Gamma^{\,\mathrm{p}}_\mathrm{in}$, respectively. Furthermore, the maximum volume is limited to $40\%$ of the fixed design domain $D$, and the weighting parameter $\omega_p$ is set to $0.1$.

\begin{figure}[H]
    \centering
    \begin{minipage}[t]{0.46\linewidth}
        \centering
        \includegraphics[width=\linewidth]{ 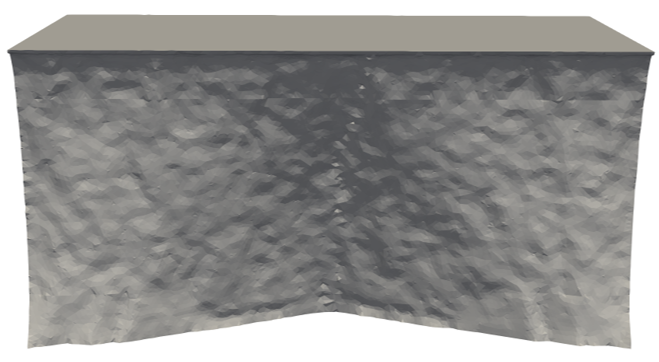}
        \subcaption{Front view without the penetrating condition}
    \end{minipage}
    \hspace{0.03\linewidth}
    \begin{minipage}[t]{0.46\linewidth}
        \centering
        \includegraphics[width=\linewidth]{ 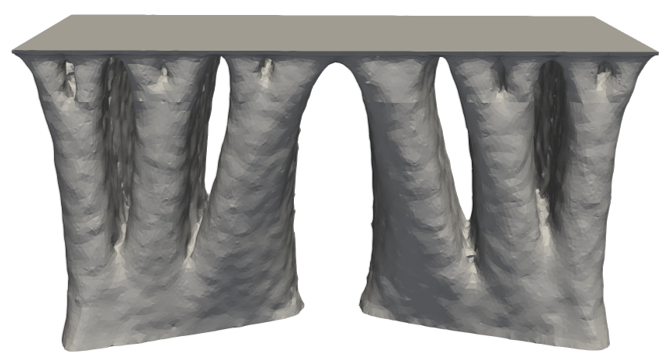}
        \subcaption{Front view with the penetrating condition}
    \end{minipage}
    \vspace{0.02\linewidth}
    \begin{minipage}[t]{0.46\linewidth}
        \centering
        \includegraphics[width=\linewidth]{ 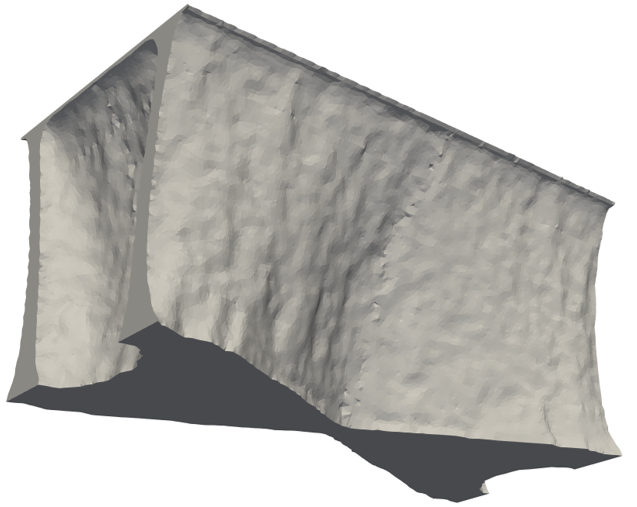}
        \subcaption{Oblique view without the penetrating condition}
    \end{minipage}
    \hspace{0.02\linewidth}
    \begin{minipage}[t]{0.46\linewidth}
        \centering
        \includegraphics[width=\linewidth]{ 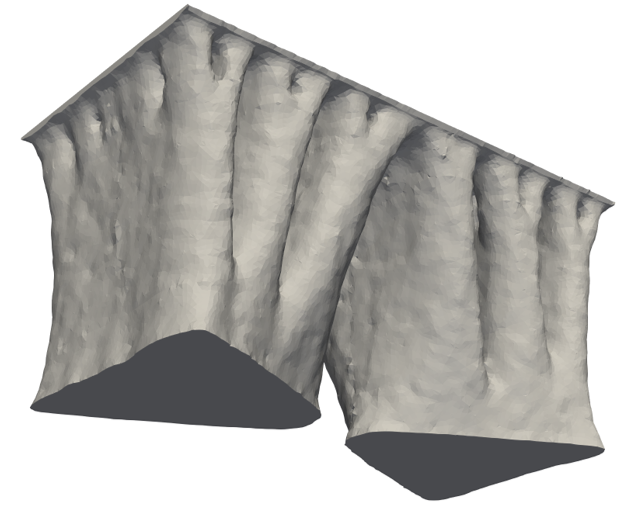}
        \subcaption{Oblique view with the penetrating condition}
    \end{minipage}
    \caption{Comparison of optimized 3D configurations without and with the penetrating condition.}
    \label{fig:3Popt}
\end{figure}

Figure \ref{fig:3Popt} compares the optimized structures without and with the penetrating condition.
As shown in (a) and (c), in the case without the penetrating condition, the optimized structure forms a shielding wall that obstructs the area between the boundaries $\Gamma^{\,\mathrm{p}}_\mathrm{out}$ and $\Gamma^{\,\mathrm{p}}_\mathrm{in}$.
However, as shown in (b) and (d), in the case with the penetrating condition, the optimized structure is penetrated.
This result demonstrates that the proposed method is useful for penetrating even in three-dimensional problems.

\begin{figure}[H]
    \begin{center}
        \begin{minipage}[t]{0.32\columnwidth}
            \centering
            \includegraphics[width=0.95\columnwidth]{ 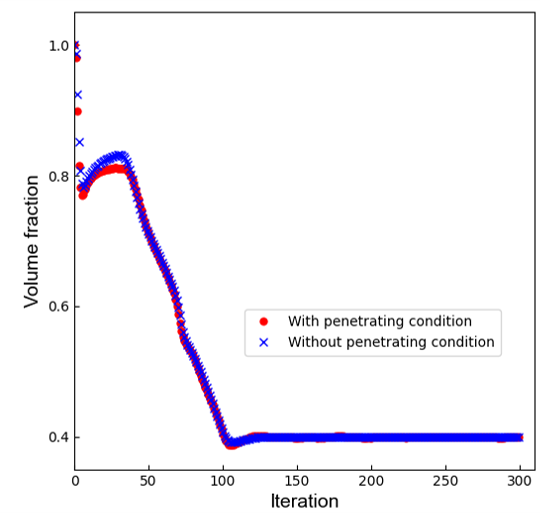}
            \subcaption{Change in volume fraction}
        \end{minipage}
        \begin{minipage}[t]{0.32\columnwidth}
            \centering
            \includegraphics[width=0.95\columnwidth]{ 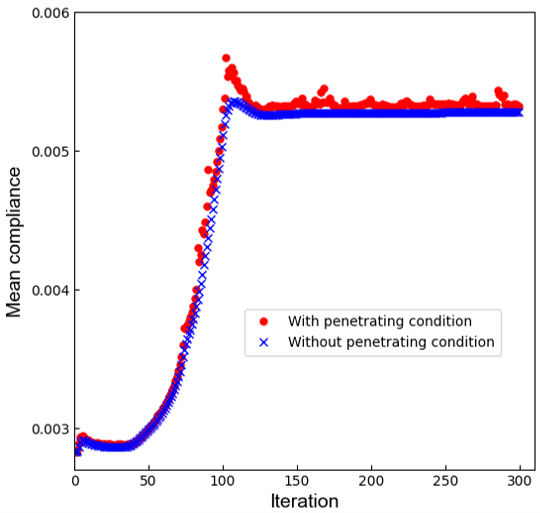}
            \subcaption{Change in mean compliance $J_u$}
        \end{minipage}
        \begin{minipage}[t]{0.32\columnwidth}
            \centering
            \includegraphics[width=0.95\columnwidth]{ 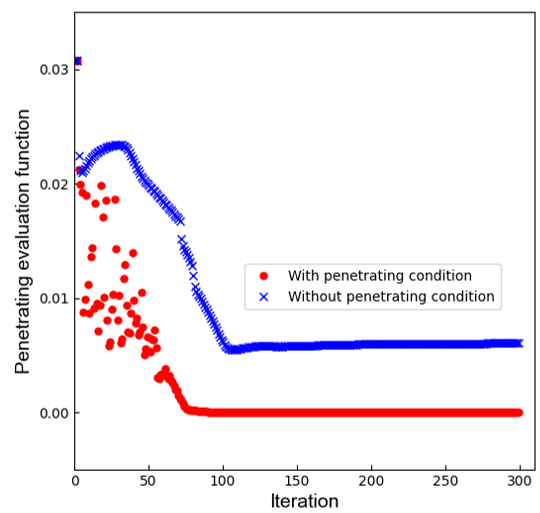}
            \subcaption{Change in penetrating evaluation function $J_p$}
        \end{minipage}
        \caption{Convergence histories comparing optimization without and with the penetrating condition.}
    \label{fig:graph3P}
    \end{center}
\end{figure}

Next, Figure \ref{fig:graph3P} shows the convergence histories of these optimization examples.
As shown in (a), both optimizations are implemented to converge the volume within approximately 100 iterations out of the total 300 iterations.
(b) illustrates the change in mean compliance for each case, and the value increases slightly by considering the shielding condition, but there is little difference between them.
In contrast, (c) shows the change in the values of the penetrating evaluation function, and there is significant difference between the cases.
Specifically, without the penetrating condition, the value remains greater than zero since the structure is not penetrated.
However, with the penetrating condition, it works to decrease the value. As a result, the structure becomes penetrated.
\subsection{2D Example Under Shielding and Penetrating Conditions} \label{sec:2D-SP}
\begin{figure}[H]
    \begin{center}
        \begin{minipage}[t]{0.49\columnwidth}
            \centering
            \includegraphics[width=0.95\columnwidth]{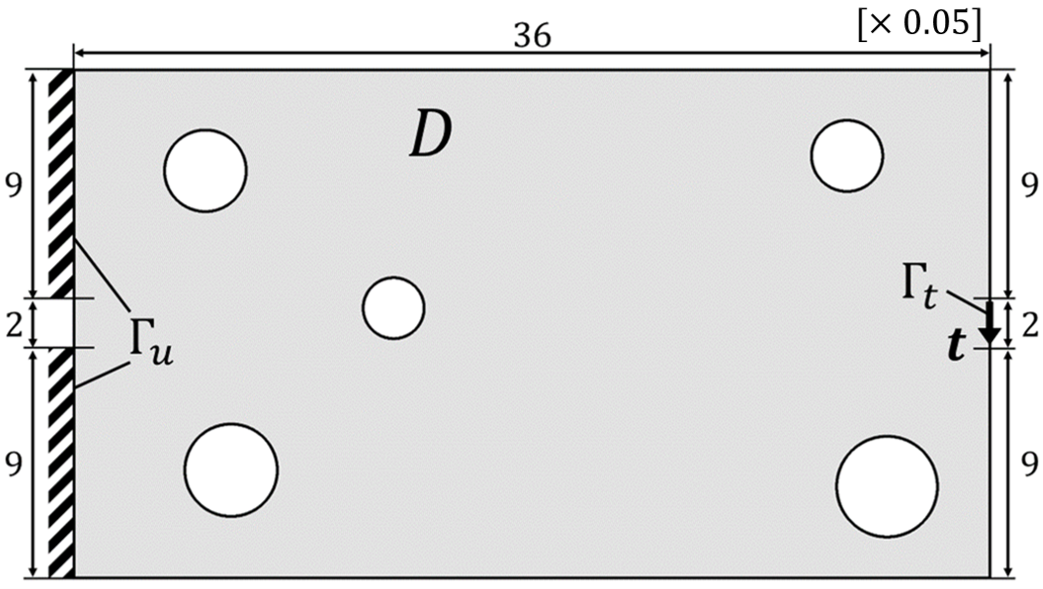}
            \subcaption{Boundary conditions for the displacement field $\bm{u}$}
        \end{minipage}
        \begin{minipage}[t]{0.49\columnwidth}
            \centering
            \includegraphics[width=0.95\columnwidth]{ 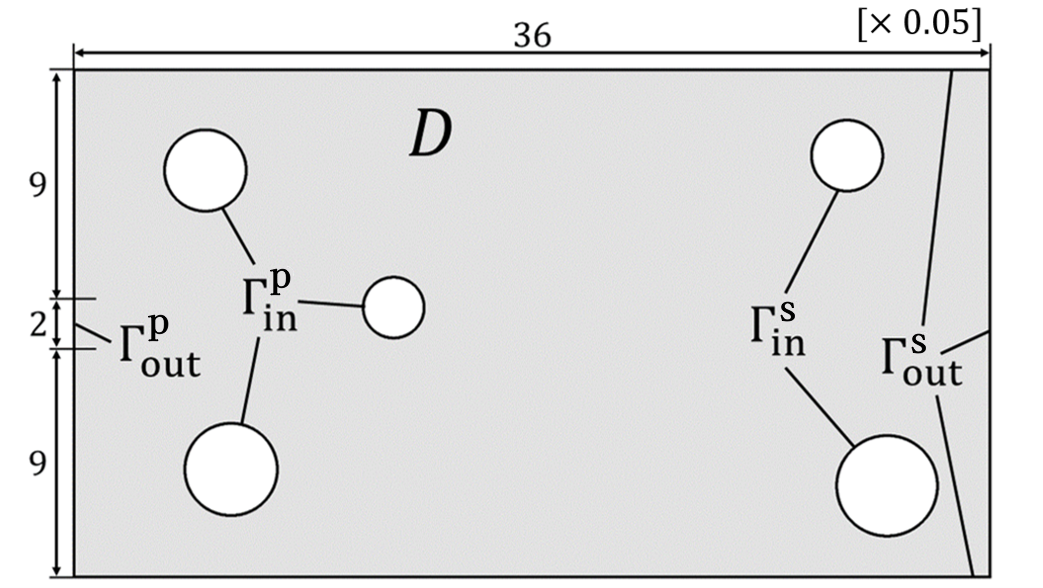}
            \subcaption{Boundary conditions for the fictitious physical fields $s$ and $p$}
        \end{minipage}
        \caption{Problem setting for mean compliance minimization considering the shielding and penetrating conditions.}
	\label{fig:2SPset}
    \end{center}
\end{figure}

Here we present the two-dimensional cantilever beam optimization as an example considering both shielding and penetrating conditions.
The problem setup is as shown in Figure \ref{fig:2SPset}, where (a) and (b) respectively indicate conditions for the displacement field and the fictitious physical fields.
The displacement is fixed at the upper and lower boundaries $\Gamma_u$ of the left edge, and the traction $\bm{t}$ is applied at the boundary $\Gamma_t$ on the middle of the right edge. As fictitious heat sources, the boundary $\Gamma^{\,\mathrm{s}}_\mathrm{out}$ is set on the upper, right, and lower edges, and the boundary $\Gamma^{\,\mathrm{p}}_\mathrm{out}$ is set on the middle of the left edge.
Additionally, for the five circular boundaries within the fixed design domain $D$, the right two are designated as the boundaries $\Gamma^{\,\mathrm{s}}_\mathrm{in}$, while the remaining three are designated as the boundaries $\Gamma^{\,\mathrm{p}}_\mathrm{in}$.
Furthermore, the maximum volume is limited to $50\%$ of the fixed design domain $D$, and the weighting parameters $\omega_s$ and $\omega_p$ are each set to  $0.05$.

\begin{figure}[H]
    \begin{center}
        \begin{minipage}[t]{0.49\columnwidth}
            \centering
            \includegraphics[width=0.95\columnwidth]{ 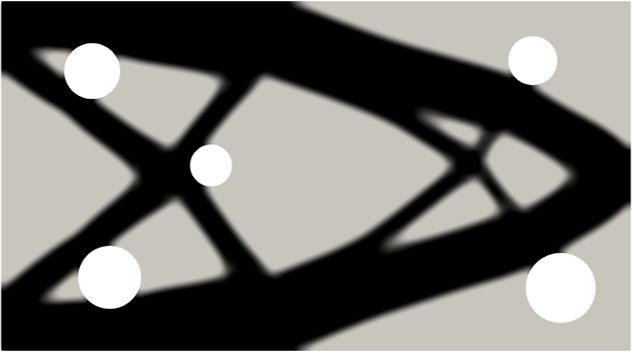}
            \subcaption{Configuration of optimized $\phi$ without shielding and penetrating conditions}
        \end{minipage}
        \begin{minipage}[t]{0.49\columnwidth}
            \centering
            \includegraphics[width=0.95\columnwidth]{ 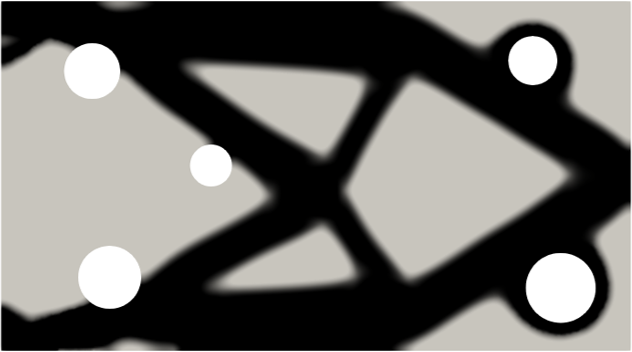}
            \subcaption{Configuration of optimized $\phi$ with shielding and penetrating conditions}
        \end{minipage}
        \caption{Comparison of 2D configurations of optimized $\phi$ without and with the shielding and penetrating conditions.}
	\label{fig:2SPopt}
    \end{center}
\end{figure}

(a) and (b) in Figure \ref{fig:2SPopt} illustrate optimized configurations without and with geometric conditions, respectively.
Firstly, focusing on the shielding feature, without the geometric conditions, the optimized structure does not shield between the boundary $\Gamma^{\,\mathrm{s}}_\mathrm{out}$ and the boundaries $\Gamma^{\,\mathrm{s}}_\mathrm{in}$. In contrast, with the geometric conditions, it shields by covering around the boundaries $\Gamma^{\,\mathrm{s}}_\mathrm{in}$.
Secondly, focusing on the penetrating feature, without the geometric conditions, the optimized structure shields between the boundary $\Gamma^{\,\mathrm{p}}_\mathrm{out}$ and the boundaries $\Gamma^{\,\mathrm{p}}_\mathrm{in}$. In contrast, with the geometric conditions, these boundaries are included within the same void domain.
These geometric differences in the optimized structures demonstrate the validity of the proposed methods, even in numerical examples that consider both shielding and penetrating conditions.
\color{black}
\subsection{3D Example Under Shielding and Penetrating Conditions} \label{sec:3D-SP}
\begin{figure}[H]
    \centering
    \begin{minipage}[t]{0.49\linewidth}
        \centering
        \includegraphics[width=\linewidth]{ 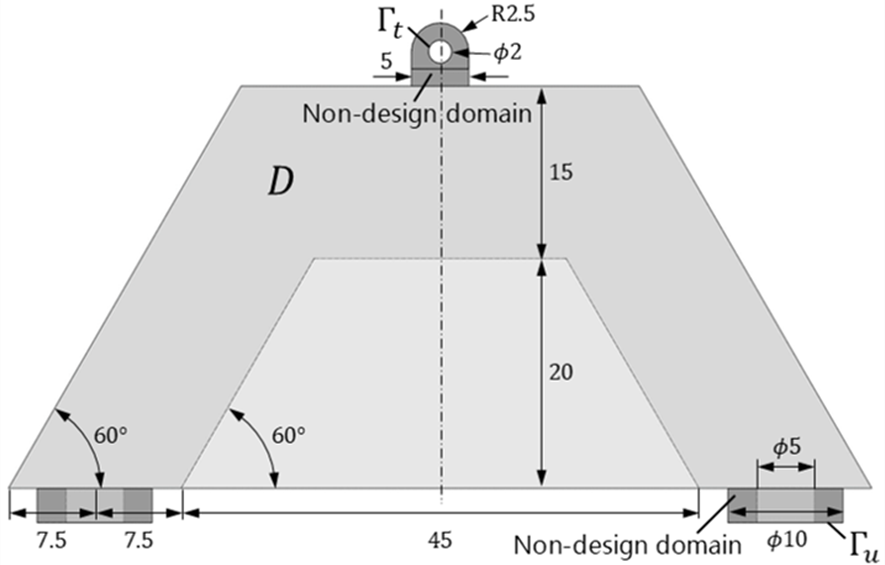}
        \subcaption{Front view of the design domain}
    \end{minipage}
    \begin{minipage}[t]{0.49\linewidth}
        \centering
        \includegraphics[width=\linewidth]{ 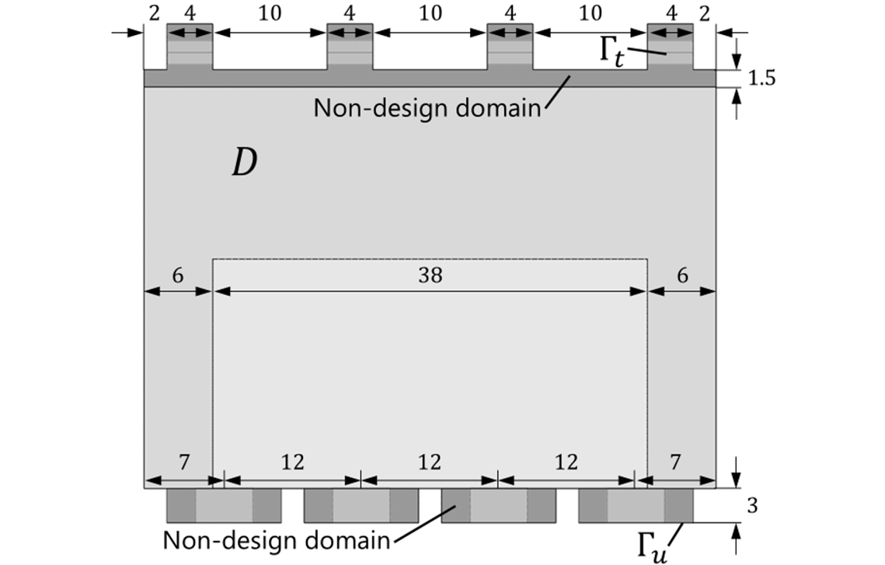}
        \subcaption{Side view of the design domain}
        \vspace{0.02\linewidth}
    \end{minipage}
    \begin{minipage}[t]{0.49\linewidth}
        \centering
        \includegraphics[width=\linewidth]{ 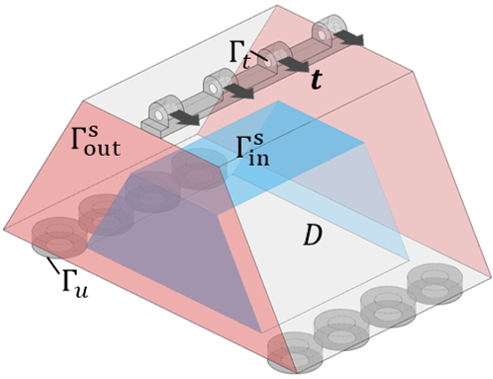}
        \subcaption{Boundary conditions for the displacement field $\bm{u}$ and the fictitious physical field $s$}
    \end{minipage}
    \begin{minipage}[t]{0.49\linewidth}
        \centering
        \includegraphics[width=\linewidth]{ 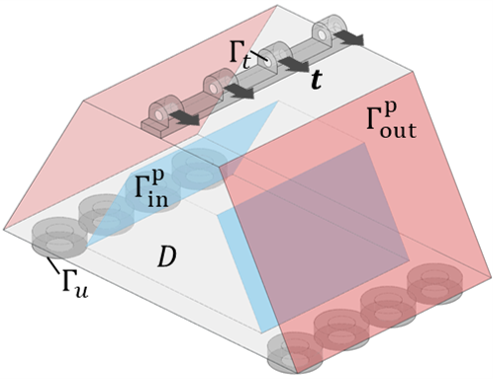}
        \subcaption{Boundary conditions for the displacement field $\bm{u}$ and the fictitious physical field $p$}
    \end{minipage}
    \caption{Problem setting for 3D mean compliance minimization considering the shielding and penetrating conditions.}
    \label{fig:3SPset}
\end{figure}

The final example is the three-dimensional bracket optimization considering both shielding and penetrating conditions.
The problem setting is shown in Figure \ref{fig:3SPset}.
Here, (a) and (b) indicate the design domain, and (c) and (d) represent the boundary conditions for the fictitious physical fields $s$ and $p$, respectively.
The design domain \(D\) is hollowed out on the underside. Fixing components are placed on the bottom surface of the design domain, and load-bearing components are placed on the upper surface. These components are distinguished as non-design domains.
The bottom surfaces of the fixing components are set to \(\Gamma_u\), and the inner surfaces of the load-bearing components are set to \(\Gamma_t\).
The traction \(\bm{t}\) pulls the boundary \(\Gamma_t\) to the right from the front view.
Additionally, regarding the fictitious physical fields, the boundary condition for $s$ is set to promote shielding in the front direction, while the boundary condition for $p$ is set to promote penetrating in the side direction.
Furthermore, the maximum volume is limited to $40\%$ of the total domain. The weighting parameters $\omega_s$ and $\omega_p$ are set to $0.15$ and $0.05$, respectively, and the regularization parameter $\tau$ is set to $5.0 \times 10^{-2}$.

\begin{figure}[H]
    \centering
    \begin{minipage}[t]{0.45\linewidth}
        \centering
        \includegraphics[width=\linewidth]{ 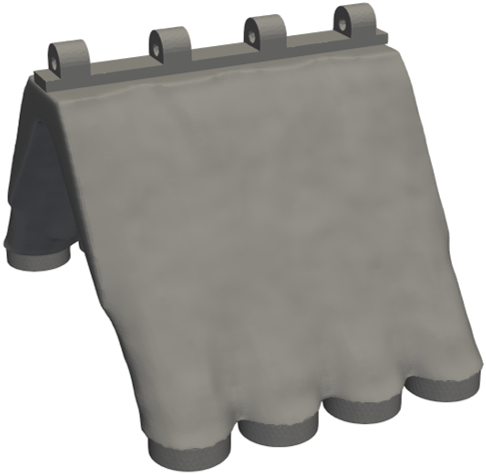}
        \subcaption{Oblique view without shielding and penetrating conditions}
    \end{minipage}
    \begin{minipage}[t]{0.45\linewidth}
        \centering
        \includegraphics[width=\linewidth]{ 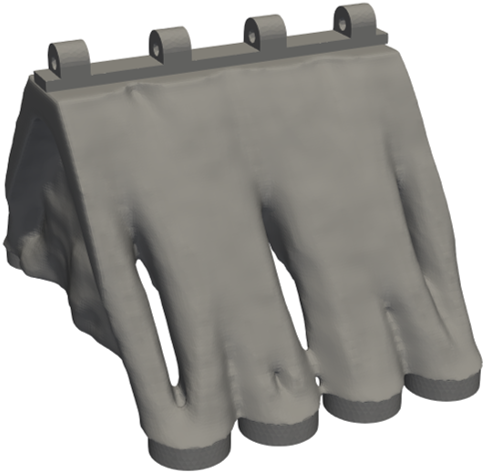}
        \subcaption{Oblique view with shielding and penetrating conditions}
        \vspace{0.02\linewidth}
    \end{minipage}
    \begin{minipage}[t]{0.49\linewidth}
        \centering
        \includegraphics[width=\linewidth]{ 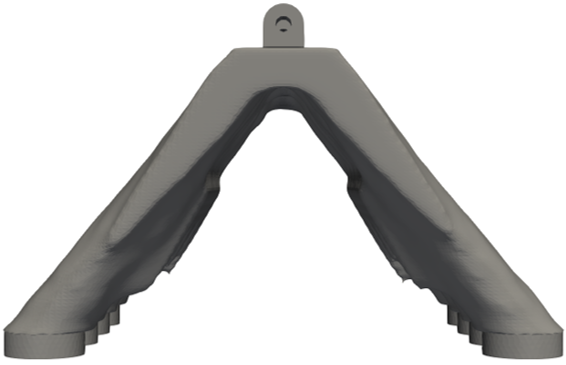}
        \subcaption{Front view without shielding and penetrating conditions}
    \end{minipage}
    \begin{minipage}[t]{0.49\linewidth}
        \centering
        \includegraphics[width=\linewidth]{ 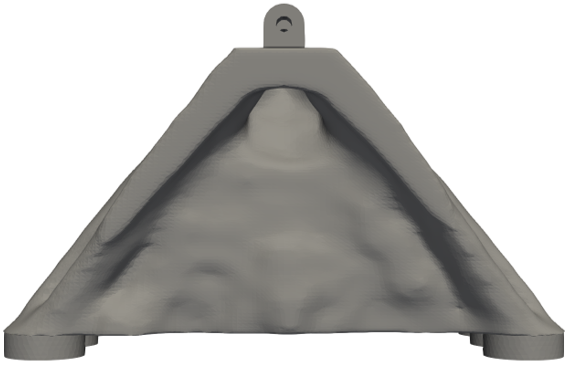}
        \subcaption{Front view with shielding and penetrating conditions}
        \vspace{0.02\linewidth}
    \end{minipage}
    \begin{minipage}[t]{0.49\linewidth}
        \centering
        \includegraphics[width=\linewidth]{ 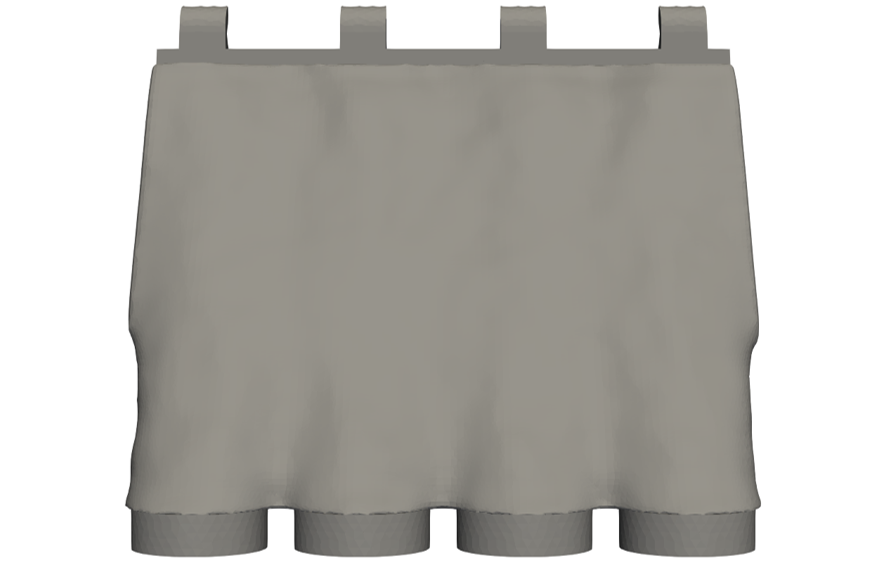}
        \subcaption{Side view without shielding and penetrating conditions}
    \end{minipage}
    \begin{minipage}[t]{0.49\linewidth}
        \centering
        \includegraphics[width=\linewidth]{ 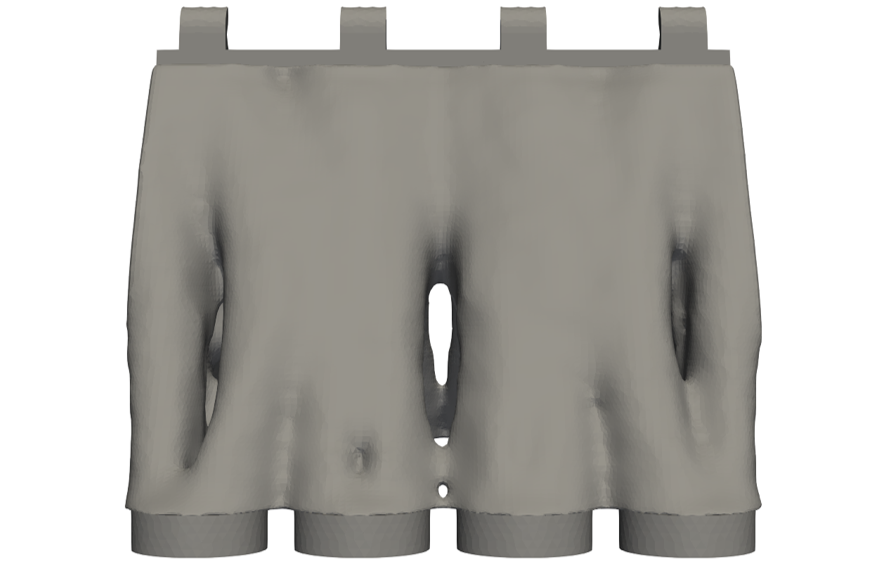}
        \subcaption{Side view with shielding and penetrating conditions}
    \end{minipage}
    \caption{\color{red}Comparison of optimized 3D configurations without and with the shielding and penetrating conditions.}
    \label{fig:3SPopt}
\end{figure}

(a),(c),(e) and (b),(d),(f) in Figure \ref{fig:3SPopt} illustrate optimized configurations without and with geometric conditions, respectively.
Regarding the shielding feature, without the geometric conditions, the optimized structure does not shield in the front direction, whereas with the geometric conditions, shielding walls are formed between boundaries $\Gamma^{\,\mathrm{s}}_\mathrm{out}$ and  $\Gamma^{\,\mathrm{s}}_\mathrm{in}$.
On the other hand, regarding the penetrating feature, without the geometric conditions, the optimized structure does not penetrate in the side direction, whereas with the geometric conditions, penetrating holes are formed between boundaries $\Gamma^{\,\mathrm{p}}_\mathrm{out}$ and  $\Gamma^{\,\mathrm{p}}_\mathrm{in}$.
These geometric differences demonstrate the validity of the proposed methods, even in three-dimensional examples that consider both shielding and penetrating conditions.
\color{black}
\section{Conclusion}
This paper proposes topology optimization that considers shielding and penetrating conditions based on the fictitious physical model.
The key points are summarized as follows:
\begin{enumerate}
    \item The framework of level set based topology optimization is described, and the fictitious physical model for considering geometric features is discussed.
    \item The fictitious physical model is proposed to evaluate the shielding and penetrating features of structures.
    \item Topology optimization for the minimum mean compliance problem considering shielding and penetrating conditions is formulated.
    \item The effects of the parameters used in the proposed method on geometric features are verified through numerical examples.
    \item It is confirmed from 2D and 3D optimization examples that the proposed method can control target geometric features.
\end{enumerate}
\section*{Acknowledgements}
This work was supported in part by JSPS KAKENHI Grant Number JP23H03800.

\bibliographystyle{elsarticle-num-names} 
\bibliography{mybibfile}
\end{document}